\documentclass{article}
\usepackage{authblk}
\title{{A comparison study of deep Galerkin method and deep Ritz method for elliptic problems with different boundary conditions}}
\author[ab]{Jingrun Chen \thanks{jingrunchen@suda.edu.cn}}
\author[ab]{Rui Du \thanks{durui@suda.edu.cn}}
\author[a]{Keke Wu \thanks{wukekever@gmail.com}}
\affil[a]{School of Mathematical Sciences, Soochow University, Suzhou, 215006, China}
\affil[b]{Mathematical Center for Interdisciplinary Research, Suzhou, 215006, China}

\date{\today}  
\usepackage{amsmath,amsthm,amssymb}
\usepackage{graphicx}
\usepackage{subfigure}
\usepackage{float}
\usepackage[table,xcdraw]{xcolor}
\usepackage{multirow}
\usepackage{booktabs}
\usepackage{cite}
\usepackage{epstopdf}
\usepackage{hyperref} 

\begin{document}
	\maketitle
	\markboth{A comparison study of deep Galerkin method and deep Ritz method for elliptic problems with different boundary conditions}
	{A comparison study of deep Galerkin mehtod and deep Ritz method}
	
	\begin{abstract}
		Recent years have witnessed growing interests in solving partial differential equations by deep neural networks, especially in the high-dimensional case. Unlike classical numerical methods, such as finite difference method and finite element method, the enforcement of boundary conditions in deep neural networks is highly nontrivial. One general strategy is to use the penalty method. {In the work, we conduct a comparison study for elliptic problems with four different boundary conditions, i.e., Dirichlet, Neumann, Robin, and periodic boundary conditions, using two representative methods: deep Galerkin method and deep Ritz method. In the former, the PDE residual is minimized in the least-squares sense while the corresponding variational problem is minimized in the latter. Therefore, it is reasonably expected that deep Galerkin method works better for smooth solutions while deep Ritz method works better for low-regularity solutions.} However, by a number of examples, we observe that deep Ritz method can outperform deep Galerkin method with a clear dependence of dimensionality even for smooth solutions and deep Galerkin method can also outperform deep Ritz method for low-regularity solutions. Besides, in some cases, when the boundary condition can be implemented in an exact manner, we find that such a strategy not only provides a better approximate solution but also facilitates the training process.
	\end{abstract}
	
	\section{Introduction}
	
	In the past decade, deep learning has achieved great success in many subjects, like computer vision, speech recognition, and natural language processing \cite{NIPS2012_4824,hinton2012deep,Goodfellow2016} due to the strong representability of deep neural networks (DNNs). Meanwhile, DNNs have also been used to solve partial differential equations (PDEs); see for example \cite{lagaris1998artificial,raissi2018hidden,weinan2017deep, berg2018unified,E2018,long2017pde,han2018solving,deepGalerkin2018,zang2020weak}. In classical numerical methods such as finite difference method \cite{leveque2007finite} and finite element method \cite{brenner2007mathematical}, the number of degrees of freedoms (dofs) grows exponentially fast as the dimension of PDE increases. One striking advantage of DNNs over classical numerical methods is that the number of dofs only grows (at most) polynomially. Therefore, DNNs are particularly suitable for solving high-dimensional PDEs. The magic underlying this is to approximate a function using the network representation of independent variables without using mesh points. Afterwards, Monte-Carlo method is used to approximate the loss (objective) function which is defined over a high-dimensional space. Some methods are based on the PDE itself \cite{raissi2018hidden,deepGalerkin2018} 
	and some other methods are based on the variational or the weak formulation {\cite{E2018,kharazmi2019variational,liao2019deep,zang2020weak}. Another successful example is the multilevel Picard approximation which is provable to overcome the curse of dimensionality for a class of semilinear parabolic equations~\cite{hutzenthaler2020proof}.} In the current work, we focus on two representative methods: deep Ritz method (DRM) proposed by E and Yu \cite{E2018} and deep Galerkin method (DGM) proposed by Sirignano and Spiliopoulos \cite{deepGalerkin2018}. {It is worth mentioning that the loss function in DGM is defined as the PDE residual in the least-squares sense. Therefore, DGM is not a Galerkin method and has no connection with Galerkin from the perspective of numerical PDEs although it is named after Galerkin.} 

	In classical numerical methods, boundary conditions can be exactly enforced for mesh points at the boundary. Typically boundary conditions include Dirichlet, Neumann, Robin, and periodic boundary conditions \cite{evans_2010}. However, it is very difficult to impose exact boundary conditions for a DNN representation. Therefore, in the loss function, it is often to add a penalty term which penalizes the difference between the DNN representation on the boundary and the exact boundary condition, typically in the sense of $L^2$ norm. Only when Dirichlet boundary condition is imposed, a novel construction of two DNN representations can be used: one for the approximation of function on the boundary and the other for the approximation of function over the domain \cite{berg2018unified}. The main purpose of the current work is to provide a comprehensive study of four boundary conditions using DRM and DGM.
	The highest derivative in the loss function in DRM is lower than that in DGM, thus it is thought that DGM works better for smooth solutions while DRM works better for low-regularity solutions. However, by a number of examples, we observe that DRM can outperform DGM with a clear dependence of dimensionality even for smooth solutions and DGM can also outperform DRM for low-regularity solutions. Besides, in some cases, when the boundary condition can be implemented in an exact manner, we find that such a strategy not only provides a better approximate solution but also facilitates the training process.
		
	The paper is organized as follows. First, a brief introduction of DGM and DRM, systematic treatment of four different boundary conditions using the penalty method, and how to use DNNs to solve PDEs are given in Section \ref{sec:method}. Numerous examples with different boundary conditions are compared in Section \ref{sec:numerics}. Conclusions are drawn in Section \ref{sec:conclusion}.
	
	\section{Methodology}\label{sec:method}
	
	The usage of a DNN to solve a PDE problem consists of three parts: the loss function, the neural network structure, and the way how the loss function is optimized over the parameter space. In what follows, we first give a brief introduction of DGM and DRM. Both methods use DNNs to approximate the PDE solution, but the main difference is the choice of loss function, which is the objective function to be optimized. Afterwards, we discuss how different boundary conditions are treated using the penalty method. We then illustrate the network structure used to approximate the PDE solution. Finally, we describe the stochastic gradient descent method which is often adopted in the optimization of loss functions.
	
	\subsection{Deep Ritz method and deep Galerkin method}
	Consider the following boundary value problem over a bounded domain $\Omega\subset\mathbb{R}^d$
	\begin{equation}
	\begin{cases}
	\mathcal{L} u(x) = f(x), & \;\text{in} \;\Omega,\\
	\Gamma u(x) = g(x), & \; \text{on} \; \partial\Omega, 
	\end{cases}
	\label{eq}
	\end{equation}
	where $d$ is the dimension, $f(x)$ and $g(x)$ are given functions, $\mathcal{L}$ is a differential operator with respect to $x$, and $\Gamma$ is a boundary operator which represents Dirichlet, Neumann, Robin, or periodic boundary condition. To proceed, we assume the well-posedness of \eqref{eq}. 
	
	The basic idea of solving a PDE using DNNs is to seek an approximate solution represented by a DNN in a certain sense \cite{hornik1989multilayer}. 
	Denote the approximate solution by $u(x;\theta)$ with $\theta$ the set of neural network parameters. Both DRM and DGM use DNNs to approximate the solution, and they only differ by the corresponding loss function. Precisely, loss functions associated to DGM and DRM in terms of $u(x;\theta)$ read as 
	\begin{equation*}
	\mathcal{J}_{\textrm{DGM}} [u(x;\theta)] = \int_{\Omega} {|\mathcal{L} u(x;\theta) - f(x)|}^2 \mathrm{d}x,
	\end{equation*}
	and
	\begin{equation*}
	\mathcal{J}_{\textrm{DRM}} [u(x;\theta)] = \int_{\Omega} \left(W(u(x;\theta)) - f(x)u(x;\theta)\right) \mathrm{d}x,
	\end{equation*}
	respectively. DGM aims to minimize the imbalance when the approximate DNN solution is substituted into $\mathcal{L} u(x) = f(x)$ of \eqref{eq} in the least-squares sense \cite{deepGalerkin2018}. DRM works in a variational sense that the variation of $\mathcal{J}_{\textrm{DRM}} [u(x;\theta)]$ with respect to $u(x;\theta)$ yields the associated Euler-Lagrange equation $\mathcal{L} u(x) = f(x)$ \cite{E2018, liao2019deep}. 
	
	The inclusion of boundary conditions is done by adding a penalty term
	\begin{equation*}
	\mathcal{B} [u(x;\theta)] = \int_{\partial\Omega} {| \Gamma u(x;\theta) - g(x) |}^2 \mathrm{d} s,
	\end{equation*} 
	and respectively, the total loss functions $\mathcal{I} [u(x;\theta)]$ for DGM and DRM are
	\begin{equation}\label{eqn:dgm}
	\mathcal{I}_{\textrm{DGM}} [u(x;\theta)] = \mathcal{J}_{\textrm{DGM}} [u(x;\theta)]  + \lambda \mathcal{B} [u(x;\theta)],
	\end{equation}
	and
	\begin{equation}\label{eqn:drm}
	\mathcal{I}_{\textrm{DRM}} [u(x;\theta)] = \mathcal{J}_{\textrm{DRM}} [u(x;\theta)]  + \lambda \mathcal{B} [u(x;\theta)],
	\end{equation}
	where $\lambda$ is the penalty parameter. 
	
	The optimal approximation $u^*(x;\theta^*)$ is obtained by solving the following optimization problem:
	\begin{equation}\label{eqn:optimization}
	 u^*(x;\theta^*)  = \arg\min_{u(x;\theta) \in \mathcal{H}(\Omega)} \mathcal{I} [u(x;\theta)],
	\end{equation}
	where $\mathcal{H}(\Omega)$ is the set of admissible functions. 
	
	\subsection{Boundary conditions}
	
	To illustrate the penalty method for boundary conditions in DGM and DRM, we start with the following explicit example  over $\Omega = (0, 1)^d$ (by default)
	\begin{equation}\label{eqn:poisson}
	- \Delta u + \pi^2 u = f(x).
	\end{equation}
	The corresponding loss terms in DGM and DRM are
	\begin{equation}\label{eqn:poissondgm}
	\mathcal{J}_{\textrm{DGM}} [u(x;\theta)] = \int_{\Omega} {| - \Delta u(x;\theta) + \pi^2 u(x;\theta) - f(x) |}^2 \mathrm{d}x,
	\end{equation}
	and
	\begin{equation}\label{eqn:poissondrm}
	\mathcal{J}_{\textrm{DRM}} [u(x;\theta)] = \int_{\Omega} {\frac{1}{2} \left( {|\nabla u(x;\theta)|}^2 + \pi^2 {u(x;\theta)}^2 \right) - f(x) u(x;\theta)} \mathrm{d}x,
	\end{equation}
	respectively.
	
	For comparison, the exact solution is set to be $u(x) = \sum_{k=1}^d \cos(\pi x_k)$ which is smooth. $f(x)$ and $g(x)$ which can be calculated explicitly will be specified later .
	
	\subsubsection{Dirichlet boundary condition}
	Dirichlet boundary condition reads as
	\begin{equation*}
	u(x) = g(x), \; x\in \partial \Omega,
	\end{equation*}
	and the corresponding penalty term is
	\begin{equation}\label{eqn:poissonpenaltyd}
	\mathcal{B}_{\textrm{D}} [u(x;\theta)] = \int_{\partial\Omega} {| u(x;\theta) - g(x) |}^2 \mathrm{d}s.
	\end{equation}
	Thus, the total loss functions of DGM and DRM for Dirichlet boundary condition are  
	\begin{equation}\label{eqn:poissondgmd}
		\mathcal{I}_{\textrm{DGM}} [u(x;\theta)] = \mathcal{J}_{\textrm{DGM}} [u(x;\theta)] + \lambda \mathcal{B}_{\textrm{D}} [u(x;\theta)],
	\end{equation}
	and
	\begin{equation}\label{eqn:poissondrmd}
		\mathcal{I}_{\textrm{DRM}} [u(x;\theta)] = \mathcal{J}_{\textrm{DRM}} [u(x;\theta)] + \lambda \mathcal{B}_{\textrm{D}} [u(x;\theta)],
	\end{equation}
	respectively.
	
	\subsubsection{Neumann boundary condition}
	
	Neumann boundary condition reads as
	\begin{equation*}
	{\partial u} / {\partial { n}} = g(x), \; x\in\partial \Omega, 
	\end{equation*}
	where $ {\partial u} / {\partial { n}} := \left({\partial u} / {\partial x_1},\cdots,{\partial u} / {\partial x_d}\right) \cdot n(x)$ 
	and $n(x)$ is the unit outer normal vector along $\partial \Omega$. The corresponding penalty term is
	\begin{equation}\label{eqn:poissonpenaltyn}
		\mathcal{B}_{\textrm{N}} [u(x;\theta)] = \int_{\partial\Omega} {| {\partial u(x;\theta)} / {\partial n} - g(x) |}^2 \mathrm{d}s
	\end{equation}
	Thus, the total loss functions of DGM and DRM for Neumann boundary condition are  
	\begin{equation}\label{eqn:poissondgmn}
		\mathcal{I}_{\textrm{DGM}} [u(x;\theta)] = \mathcal{J}_{\textrm{DGM}} [u(x;\theta)] + \lambda \mathcal{B}_{\textrm{N}} [u(x;\theta)],
	\end{equation}
	and
	\begin{equation}\label{eqn:poissondrmn}
		\mathcal{I}_{\textrm{DRM}} [u(x;\theta)] = \mathcal{J}_{\textrm{DRM}} [u(x;\theta)] + \lambda \mathcal{B}_{\textrm{N}} [u(x;\theta)],
	\end{equation}
	respectively.
	
	\subsubsection{Robin boundary condition}
	
	Robin boundary condition reads as
	\begin{equation*}
		{\partial u} / {\partial { n}} + u(x) = g(x), \; x\in \partial \Omega,
	\end{equation*}
	and the corresponding penalty term is
	\begin{equation}\label{eqn:poissonpenaltyr}
		\mathcal{B}_{\textrm{R}} [u(x;\theta)] = \int_{\partial\Omega} {| {\partial u(x;\theta)} / {\partial { n}} + u(x;\theta) - g(x) |}^2 \mathrm{d}s.
	\end{equation} 
	Thus, the total loss functions of DGM and DRM for Robin boundary condition are  
	\begin{equation}\label{eqn:poissondgmr}
		\mathcal{I}_{\textrm{DGM}} [u(x;\theta)] = \mathcal{J}_{\textrm{DGM}} [u(x;\theta)] + \lambda \mathcal{B}_{\textrm{R}} [u(x;\theta)],
	\end{equation}
	and
	\begin{equation}\label{eqn:poissondrmr}
		\mathcal{I}_{\textrm{DRM}} [u(x;\theta)] = \mathcal{J}_{\textrm{DRM}} [u(x;\theta)] + \lambda \mathcal{B}_{\textrm{R}} [u(x;\theta)],
	\end{equation} 
	respectively.
	
	\subsubsection{Periodic boundary condition}\label{sec:pbcpenalty}
	
	Periodic boundary condition over the boundary of $\Omega = (-1, 1)^d$ reads as 
	\begin{equation*}
	\begin{cases}
		u(\tilde{x}_k,-1) = u(\tilde{x}_k,1), \\
		{\partial u(\tilde{x}_k,-1)}/{\partial {x_k}} = {\partial u(\tilde{x}_k,1)}/{\partial {x_k}}, \\
	\end{cases}
	\end{equation*}
	where $\tilde{x}_k =(x_1,\cdots,x_{k-1},x_{k+1},\cdots,x_d)$ for $k=1,\cdots,d$. The exact solution is still $u(x) = \sum_{k=1}^d \cos(\pi x_k)$.
	Note that the penalty term $\mathcal{B}_{\textrm{P}} [u(x;\theta)]$ in this case consists of two terms:
	\begin{align*}
		\mathcal{B}_{\textrm{P}_1} [u(x;\theta)] & = \sum_{k=1}^{d} \int_{\partial\Omega} {|u(\tilde{x}_k,-1) -u(\tilde{x}_k,1)|}^2  \mathrm{d}s, \\
		\mathcal{B}_{\textrm{P}_2} [u(x;\theta)] & = \sum_{k=1}^{d} \int_{\partial\Omega} {|{\partial u(\tilde{x}_k,-1)}/{\partial {x_k}} - {\partial u(\tilde{x}_k,1)}/{\partial {x_k}}|}^2  \mathrm{d}s.
	\end{align*}
	Thus, the corresponding loss functions of DGM and DRM for periodic boundary condition are  
	\begin{equation}\label{eqn:poissondgmp}
		\mathcal{I}_{\textrm{DGM}} [u(x;\theta)] = \mathcal{J}_{\textrm{DGM}} [u(x;\theta)] + \lambda_1 \mathcal{B}_{\textrm{P}_1} [u(x;\theta)] + \lambda_2 \mathcal{B}_{\textrm{P}_2} [u(x;\theta)],
	\end{equation}
	and
	\begin{equation}\label{eqn:poissondrmp}
	\mathcal{I}_{\textrm{DRM}} [u(x;\theta)] = \mathcal{J}_{\textrm{DRM}} [u(x;\theta)] + \lambda_1 \mathcal{B}_{\textrm{P}_1} [u(x;\theta)] + \lambda_2 \mathcal{B}_{\textrm{P}_2} [u(x;\theta)],
	\end{equation}
	where $\lambda_1$ and $\lambda_2$ are prescribed penalty parameters.
	
	\subsection{Network structure}
	
\begin{figure}[H]
		\centering
		\includegraphics[scale=0.5]{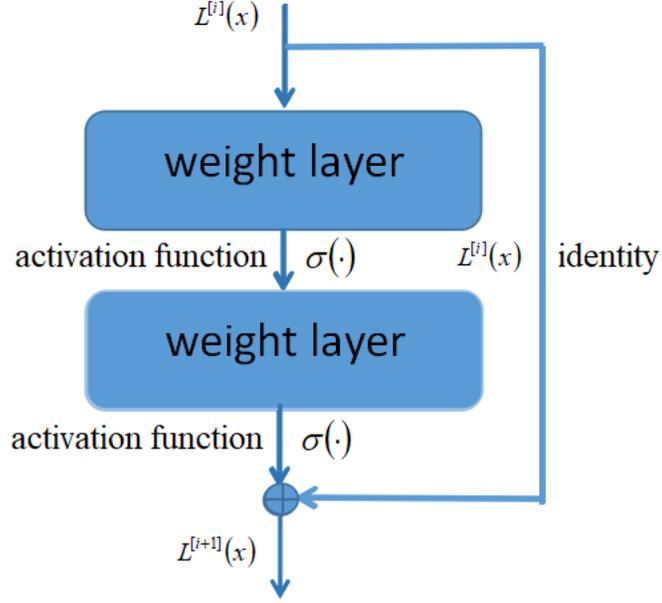}
		\caption{One residual block in the neural network structure.}
		\label{residual}
	\end{figure}

	The deep network structure employed here is similar to ResNet \cite{DBLP:journals/corr/HeZRS15}, which is built by stacking several residual blocks. Each residual block contains one input, two weight layers, and two nonlinear transformation operations (activation functions) with a skip identity connection and one output. In details, let us consider a network with $n$ residual blocks. For the $i$-th block, let $L^{[i]}(x) \in \mathbb{R}^{m}$ be the input, $W^{[i]}_1, W^{[i]}_2 \in \mathbb{R}^{m \times m}$ and $b^{[i]}_1, b^{[i]}_2 \in \mathbb{R}^{m}$ be the weight matrices and the bias vectors, $\sigma(\cdot)$ be the activation function, and $L^{[i+1]}(x)$ be the output which can be specified as
	\begin{equation}\label{eqn:resnet}
		L^{[i+1]}(x) = \sigma (W^{[i]}_2 \cdot (\sigma (W^{[i]}_1 \cdot L^{[i]}(x) + b^{[i]}_1)) + b^{[i]}_2) + L^{[i]}(x).
	\end{equation}   
	The initial input $L^{[0]}(x) = W^{[0]} \cdot x + b^{[0]}$ and the final output $L^{[n+1]}(x) = W^{[n+1]} \cdot L^{[n]}(x) + b^{[n+1]}$ with $ W^{[0]} \in \mathbb{R}^{m \times d}, b^{[0]} \in \mathbb{R}^{m \times 1}$ and $W^{[n+1]} \in \mathbb{R}^{1 \times m}, b^{[n+1]} \in \mathbb{R}$. 
	The schematic picture of one residual block is given in Figure \ref{residual}.
	{Below are activation functions used in the current work}
	\begin{align*}
		relu(x) & = \max(0,x), \\
		sigmoid(x) & = \frac{1}{1+\exp(-x)}, \\
		swish(x) & = \frac{x}{1+\exp(-x)}, \\
		\sigma(x) & = {(\sin x )}^3, \\
		swish(ax) & =  \frac{ax}{1+\exp(-ax)}.
	\end{align*}
	{The last activation function is the adaptive swish function where $a$ is an additional parameter and is optimized in the training process \cite{JAGTAP2020109136,2019arXiv190912228J}.}
	
	Overall, the DNN approximation of PDE solution can be written as
	\begin{eqnarray}\label{eqn:ddnsolution}
		u(x;\theta) =  L^{[n+1]} \circ L^{[n]} \circ \cdots \circ L^{[1]} \circ L^{[0]}(x),
	\end{eqnarray}
	where $\theta$ is the full set of all weight and bias parameters in the neural network, i.e., $\theta = \{W^{[0]}, b^{[0]}, \{W_1^{[i]}, b_1^{[i]}, W_2^{[i]},b_2^{[i]} \}_{i=1}^n, W^{[n+1]}, b^{[n+1]}\}$. The total number of parameters is $m(d+1)+(2mn+1)(m+1)$.
	
	\subsection{Stochastic gradient descent algorithm}
	
    Using DNNs to solve PDEs is now transferred to solve the optimization problem \eqref{eqn:optimization} with the loss function \eqref{eqn:dgm} or \eqref{eqn:drm} over the possible DNN representations \eqref{eqn:ddnsolution}. Even if the original PDE is linear, the DNN representation \eqref{eqn:ddnsolution} can be highly nonlinear due to the successive composition of nonlinear activation functions. On the other hand, quadrature schemes for the high-dimensional integral in \eqref{eqn:dgm} and \eqref{eqn:drm} run into the curse of dimensionality and Monte-Carlo method can overcome this issue. The stochastic gradient descent (SGD) algorithm and its variants play a key role in deep learning training. It is a first-order optimization method which naturally incorporates the idea of Monte-Carlo sampling and thus avoids the curse of dimensionality. At each iteration, SGD updates neural parameters by evaluating the gradient of the loss function only at a batch of samples as
    {
	\begin{equation}\label{eqn:sgd}
		\theta_{k+1} = \theta_{k} - \epsilon_k  \frac{1}{N}\sum_{i=1}^{N}\nabla_\theta 		
		l_i(\theta_k),
	\end{equation}}
	where $\theta_{k}$ is the parameters of neural network at the $k$-th iteration, $\epsilon_k$ is the learning rate, and {$l_i(\theta_k)$ is used to approximate the loss function using the single function value $u(x_i;\theta_k)$ multiplied by the volume or boundary measure.} $x_i$ are randomly generated with uniform distribution over $\Omega$ and $\partial\Omega$. Though better sampling strategies, such as quasi-Monte Carlo sampling \cite{chen2019quasi}, can be used, we stick to Monte-Carlo sampling \cite{Ogata1989} in the current work for the comparison purpose. 
	
	In our work, Adam optimizer is used to accelerate the training of the deep neural network \cite{ADAM}. Adam algorithm estimates first-order and second-order moments of gradient to dynamically adjust the learning rate for each parameter. The main advantage is that the learning rate at each iteration has a certain range after correction, which makes the parameter update more stable. In implementation, the global learning rate $\epsilon$ is $0.001$, the exponential decay rates of moment estimation $\rho_1, \rho_2$ are set to be $0.9$ and $0.999$, and the small constant $\delta$ used for numerical stability is set to be $10^{-8}$. 
	In addition, we use the finite difference method to approximate derivatives in the loss function. 
	
	\section{Numerical results}\label{sec:numerics}
	
	We shall use the following relative $L^2$ error to measure the approximation error
	\begin{equation*}
	\textrm{error} = \sqrt{\frac{\int_{\Omega} {\left( u^{*}(x;\theta^*) - u(x) \right)}^2 \mathrm{d} x}{\int_{\Omega} {\left( u(x) \right)}^2 \mathrm{d} x}},
	\end{equation*}
	where $u^{*}(x;\theta^*)$ is the DDN approximation of DGM or DRM and $u(x)$ is the exact solution, respectively. 
	
	\subsection{Training process and dimensional dependence for four boundary conditions}\label{sec:bcresult}
	For four different boundary conditions, we record the training process of DGM and DRM and measure the error in terms of dimensionality. For comparison purpose, the same setup is employed for different boundary conditions, but the network structure varies as the dimensionality $d$ increases. Typically, each neural network contains three to four residual blocks with several neural units in each layer. The activation function used here is $swish(x)$.
	
	{
	For Dirichlet boundary condition, there are some strategies to avoid the penalty term; see \cite{berg2018unified,sheng2020pfnn} for example. The basic idea is to employ one DNN denoted by $DNN(x;\theta)$ to approximate the PDE solution in the following trail form
	\begin{equation}\label{eqn:twostage}
	u(x;\theta) = L_D(x) DNN(x;\theta) + G(x),
	\end{equation}
	where $L_D(x)$ is the distance function to the boundary and $G(x)$ is a smooth extension of $g(x)$ over the whole domain $\Omega$. For periodic boundary condition, construction of a specific DNN can automatically satisfy the boundary condition \cite{han2020solving}. We shall return to this in Section \ref{sec:pbc}.
	For the other two boundary conditions, in principle, the strategy in \eqref{eqn:twostage} can be applied if a natural extension function $G(x)$ satisfies the boundary condition and the distance function is available. From the practical perspective, however, it is unclear how to find such an extension function $G(x)$ that satisfies the boundary condition. Therefore, we mainly focus on the penalty method for four different boundary conditions and provide results without the penalty term for both Dirichlet and periodic boundary conditions.}

\begin{figure}[H]
		\centering
		\subfigure[Dirichlet]{
			\includegraphics[width=2.0in]{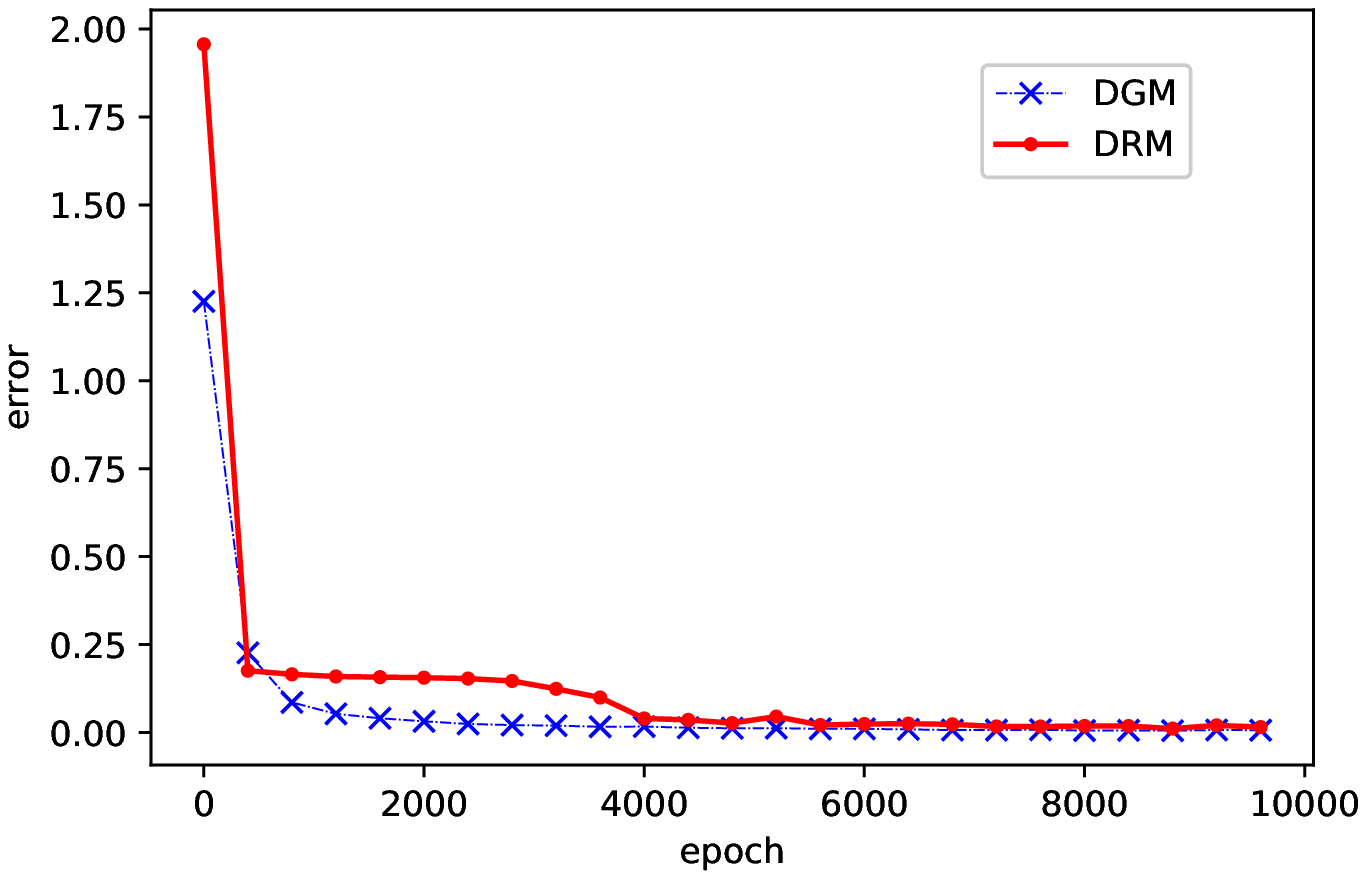}
		}
		\subfigure[Neumann]{
			\includegraphics[width=2.0in]{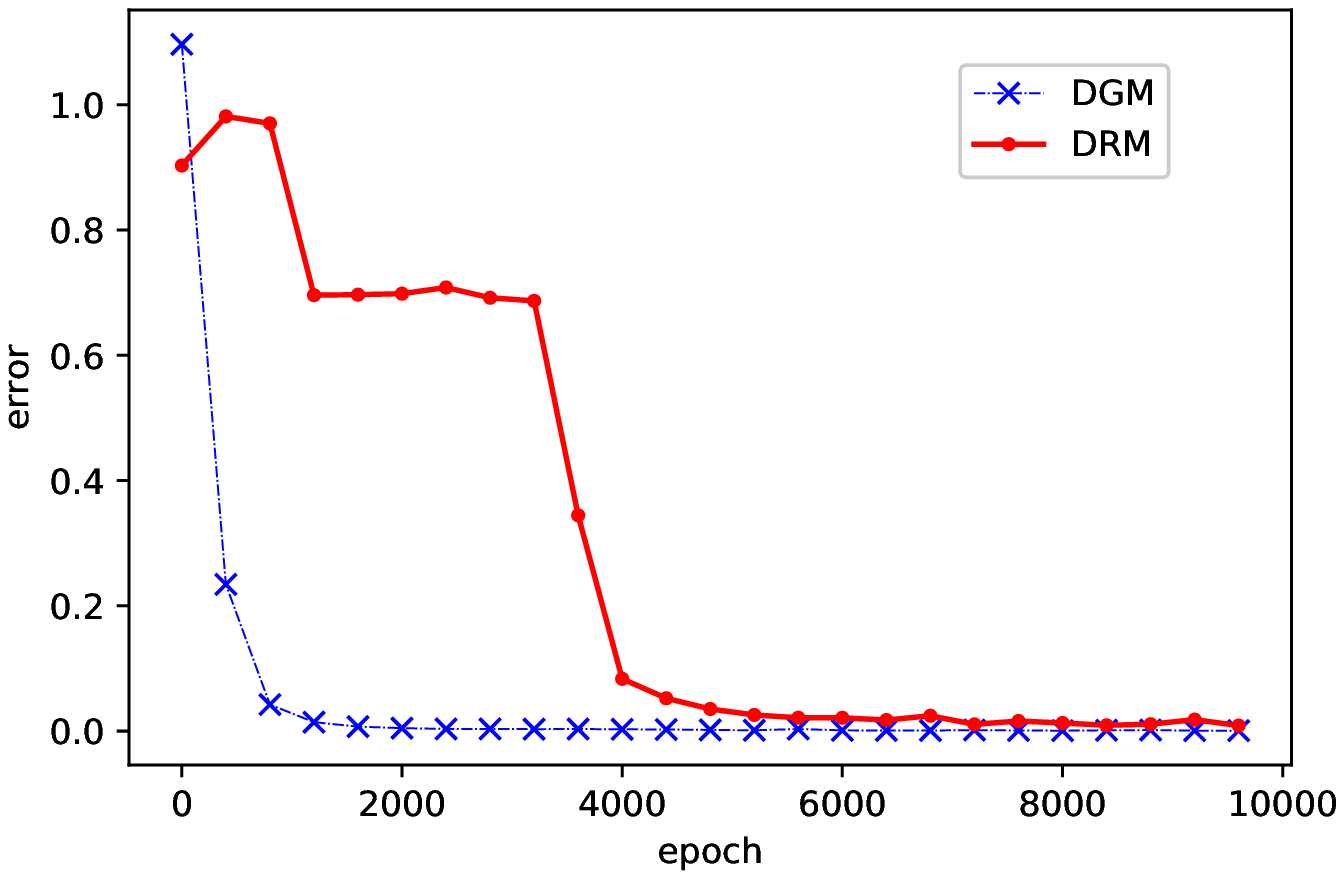}
		}
		\quad    
		\subfigure[Robin]{
			\includegraphics[width=2.0in]{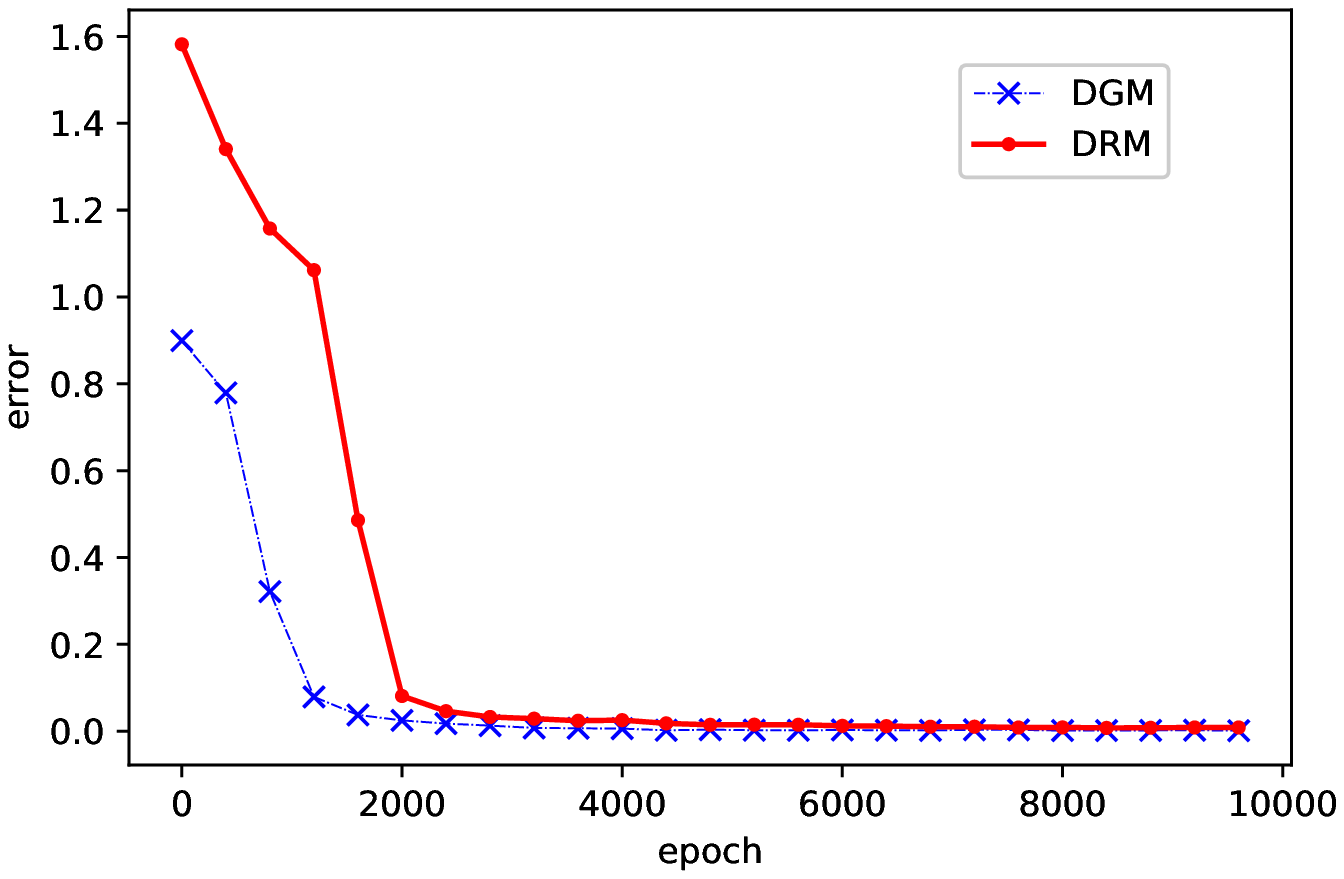}
		}
		\subfigure[Periodic]{
			\includegraphics[width=2.0in]{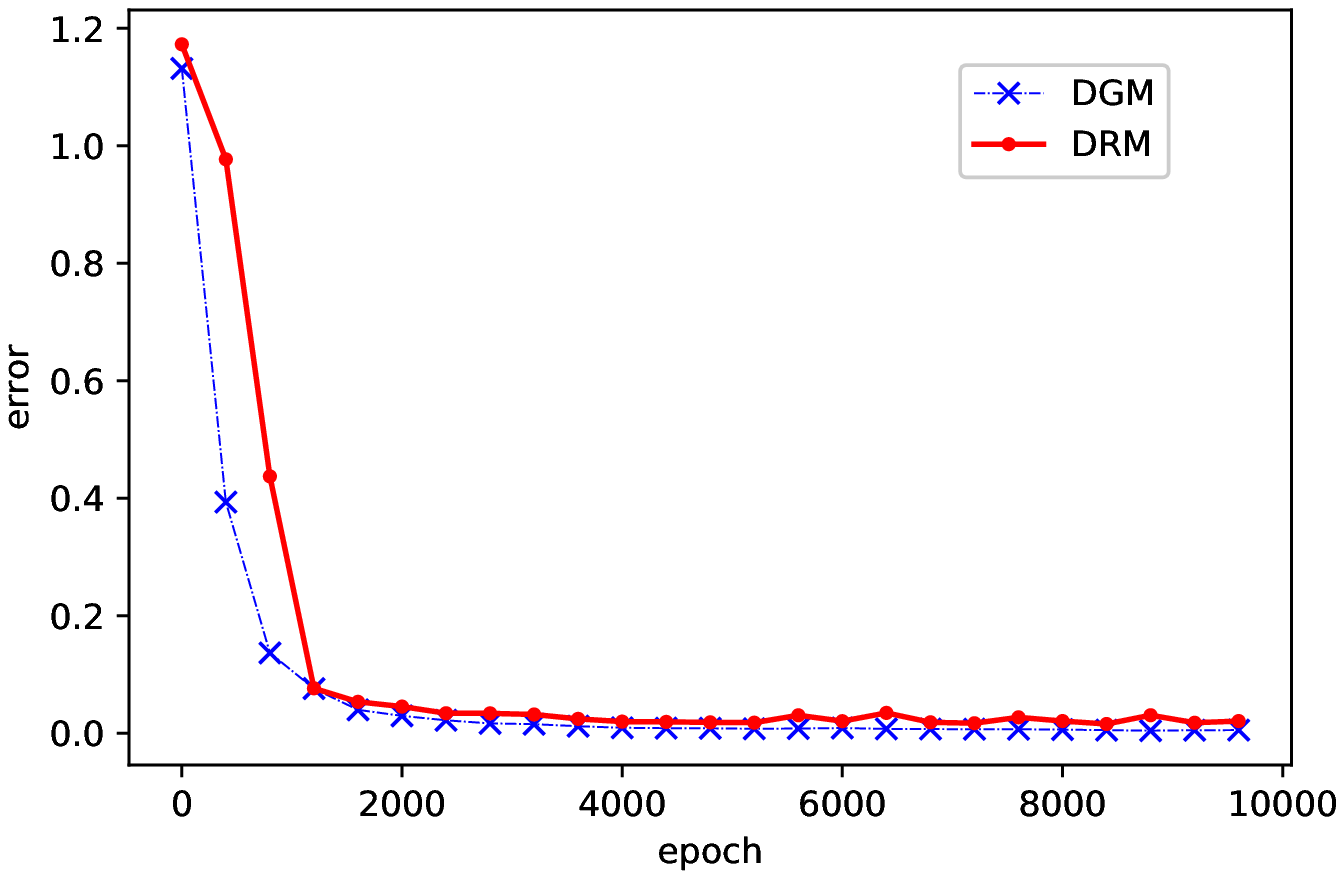}
		}
		\caption{Training processes of DGM and DRM for four boundary conditions in 2D. Each neural network contains three residual blocks with four neural units in each layer. The mini-batch size is $2000$ in the domain and $400$ on the boundary for one epoch. The penalty parameter $\lambda = 100.0$ for Dirichlet, Neumann, and Robin boundary conditions and $\lambda_1 = 10.0, \lambda_2 = 5.0$ for periodic boundary condition.}
		\label{DGM v.s. DRM in 2D}
	\end{figure}

\begin{figure}[H]
		\centering
		\subfigure[Dirichlet]{
			\includegraphics[width=2.0in]{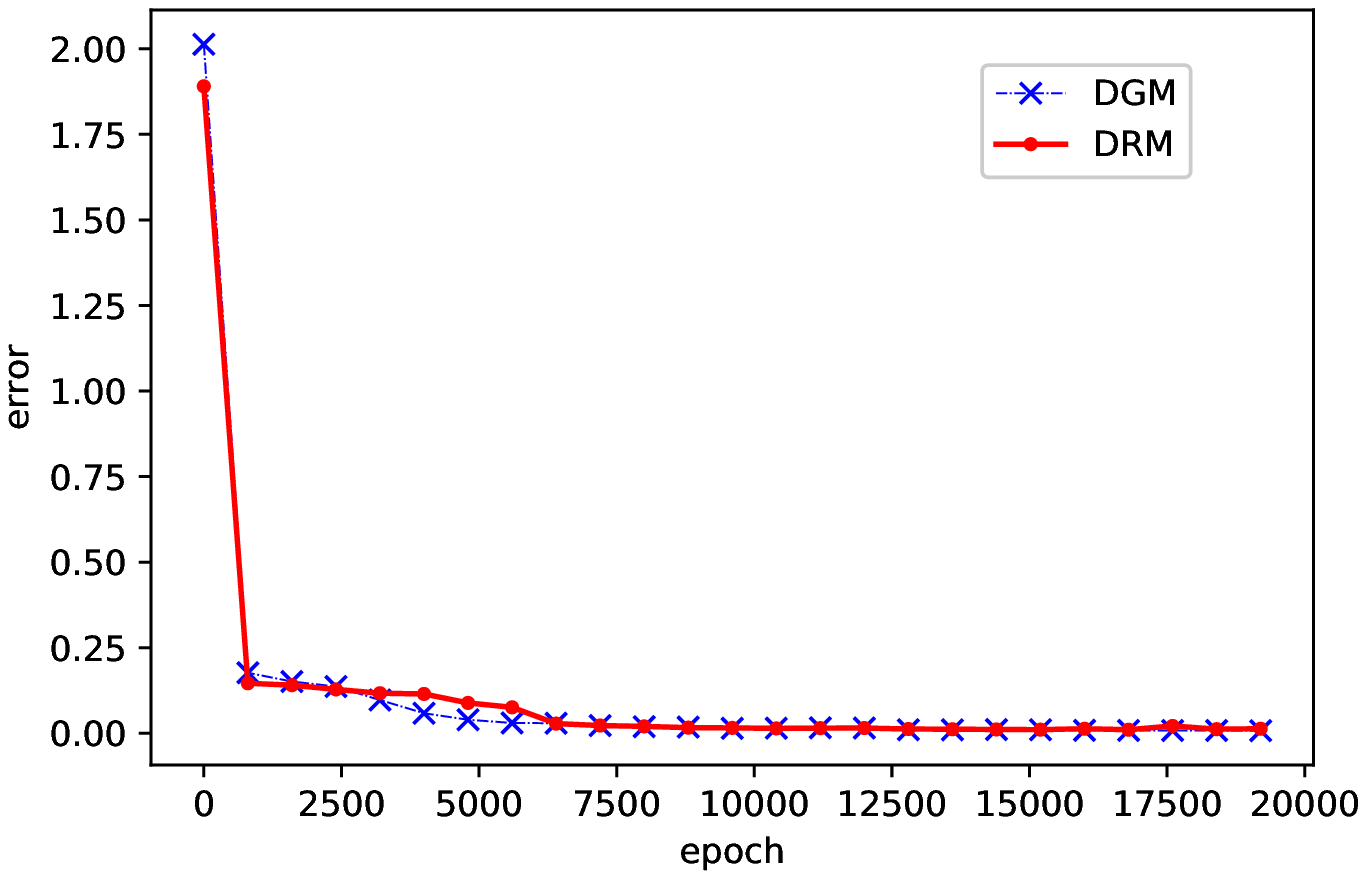}
		}
		\subfigure[Neumann]{
			\includegraphics[width=2.0in]{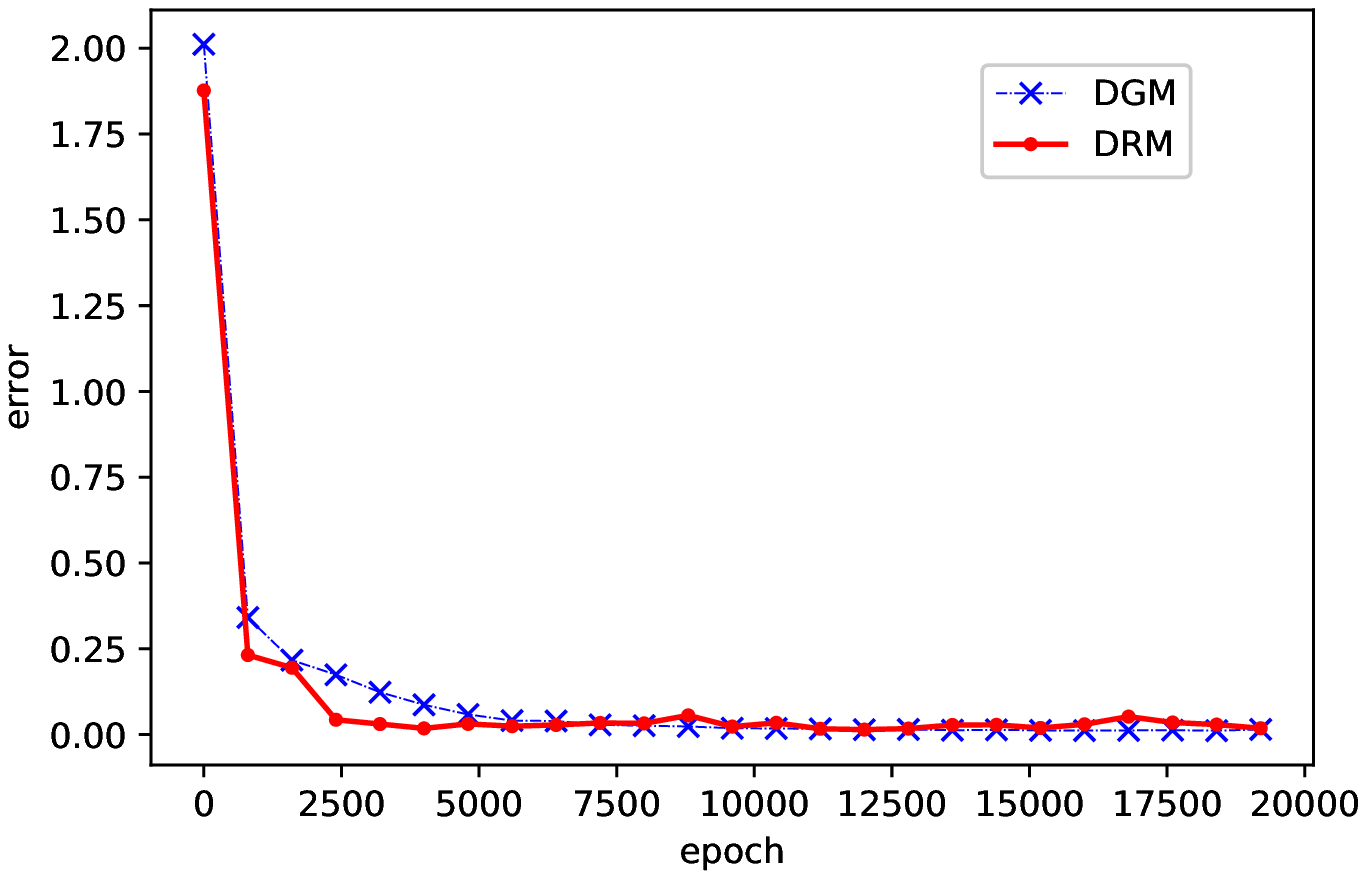}
		}
		\quad    
		\subfigure[Robin]{
			\includegraphics[width=2.0in]{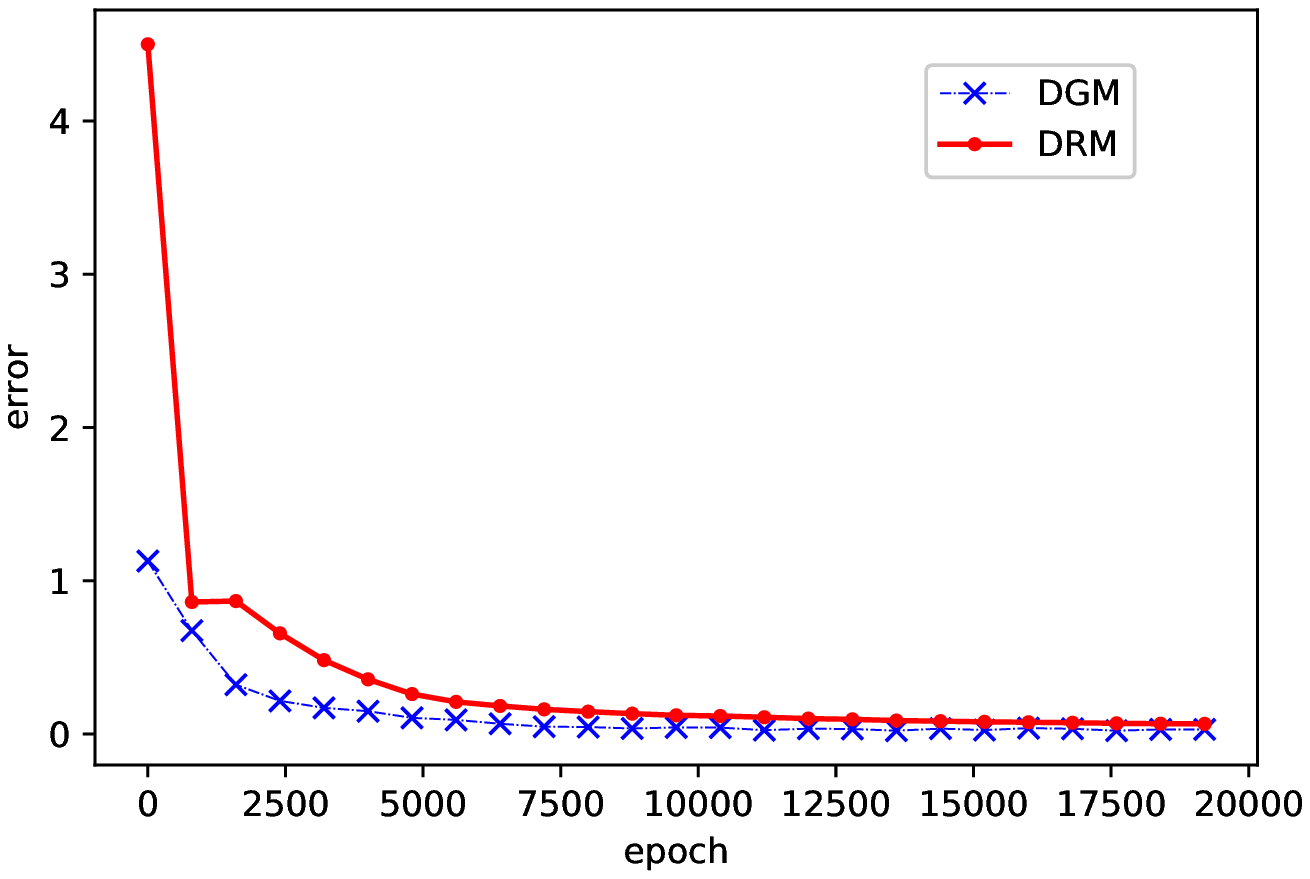}
		}
		\subfigure[Periodic]{
			\includegraphics[width=2.0in]{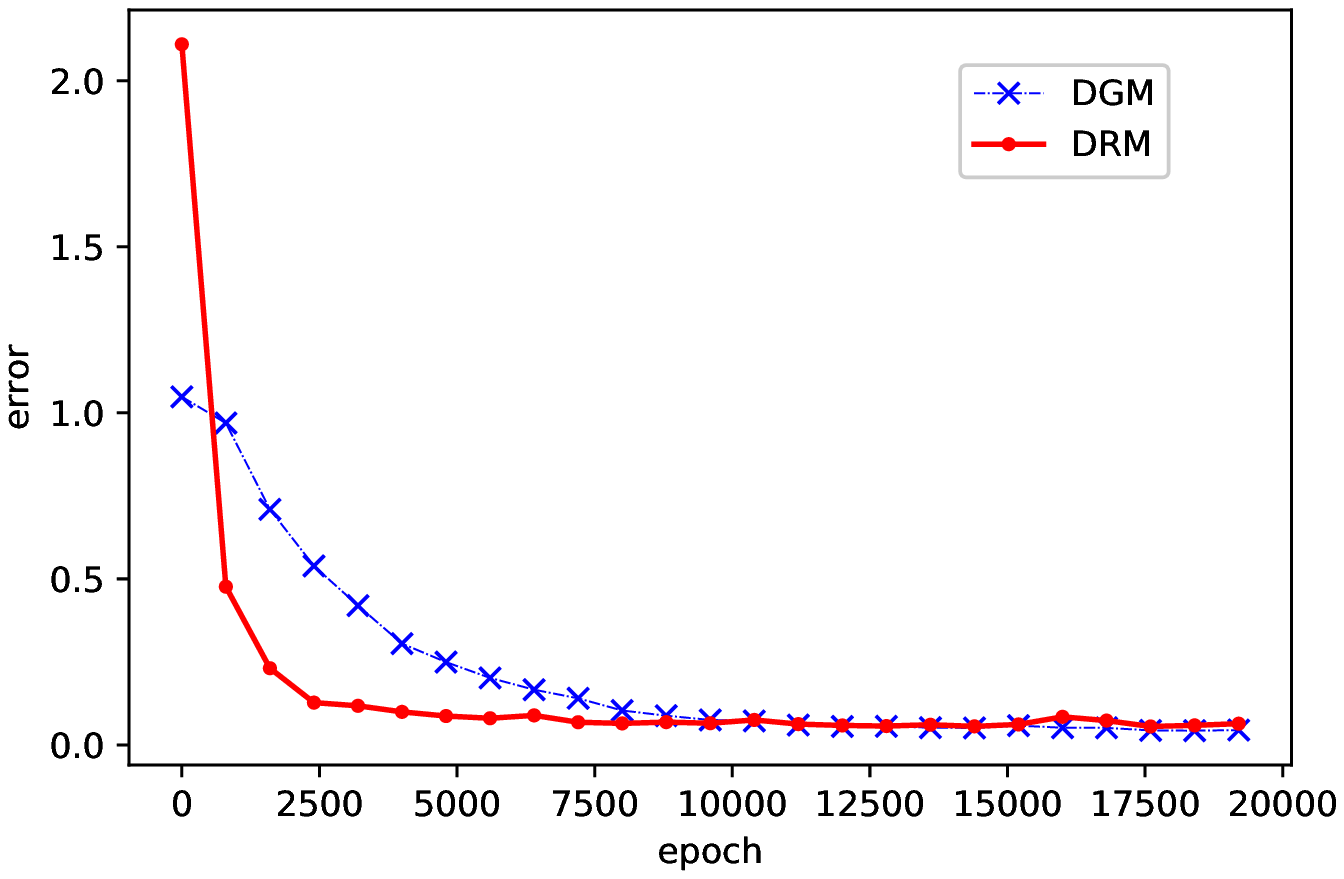}
		}
		\caption{Training processes of DGM and DRM for four boundary conditions in 4D. Each neural network contains three residual blocks with eight neural units in each layer. The mini-batch size is $2000$ in the domain, $800$ on the boundary for Dirichlet, Neumann, and Robin boundary conditions, and $8000$ on the boundary for periodic boundary condition. The penalty parameter $\lambda = 100.0$ for Dirichlet boundary condition, $\lambda = 1.0$ for Neumann boundary condition, $\lambda = 500.0$ for Robin boundary condition, and $\lambda_1 = 1.0, \lambda_2 = 0.5$ for periodic boundary condition.}
		\label{DGM v.s. DRM in 4D}
	\end{figure}

\begin{figure}[H]
		\centering
		\subfigure[Dirichlet]{
			\includegraphics[width=2.0in]{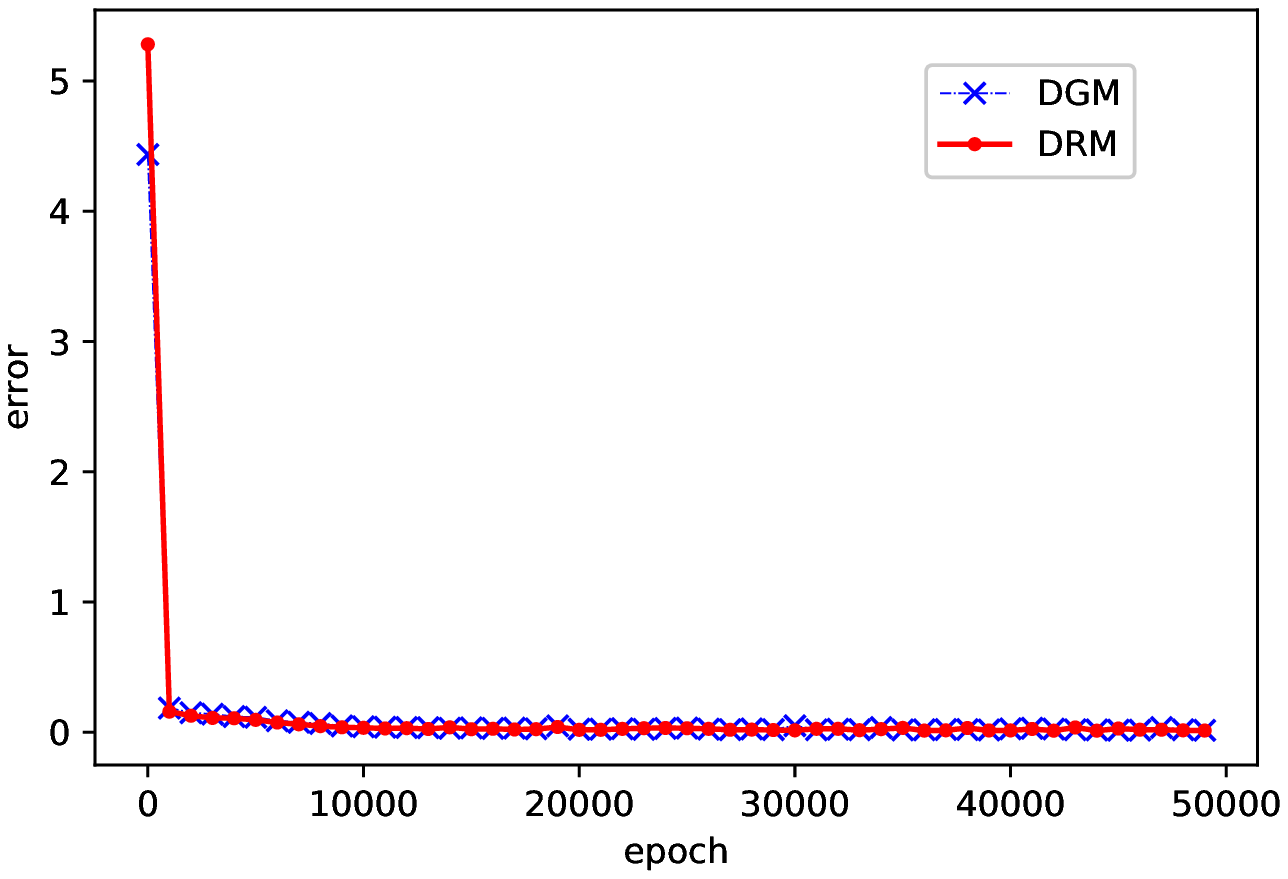}
		}
		\subfigure[Neumann]{
			\includegraphics[width=2.0in]{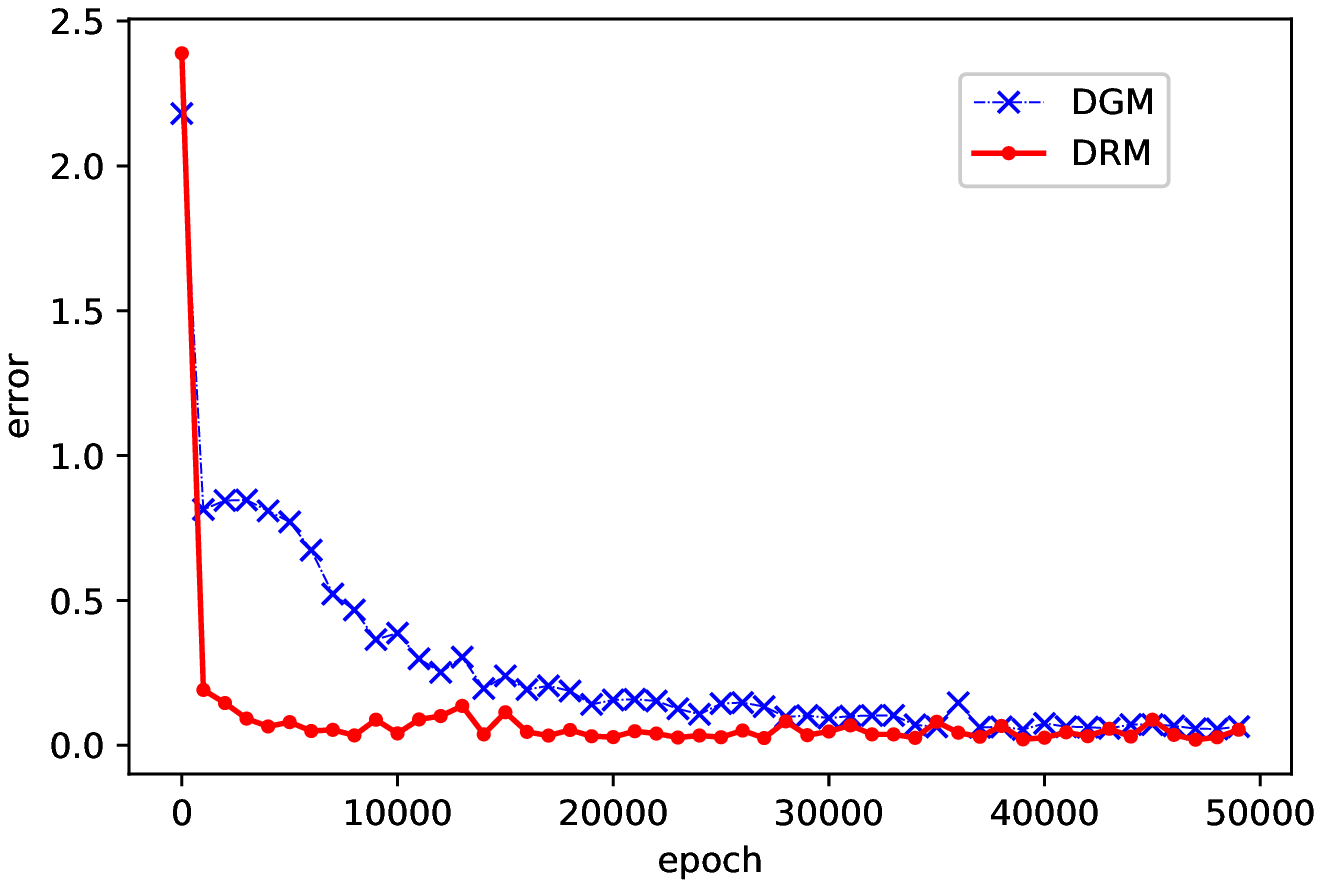}
		}
		\quad    
		\subfigure[Robin]{
			\includegraphics[width=2.0in]{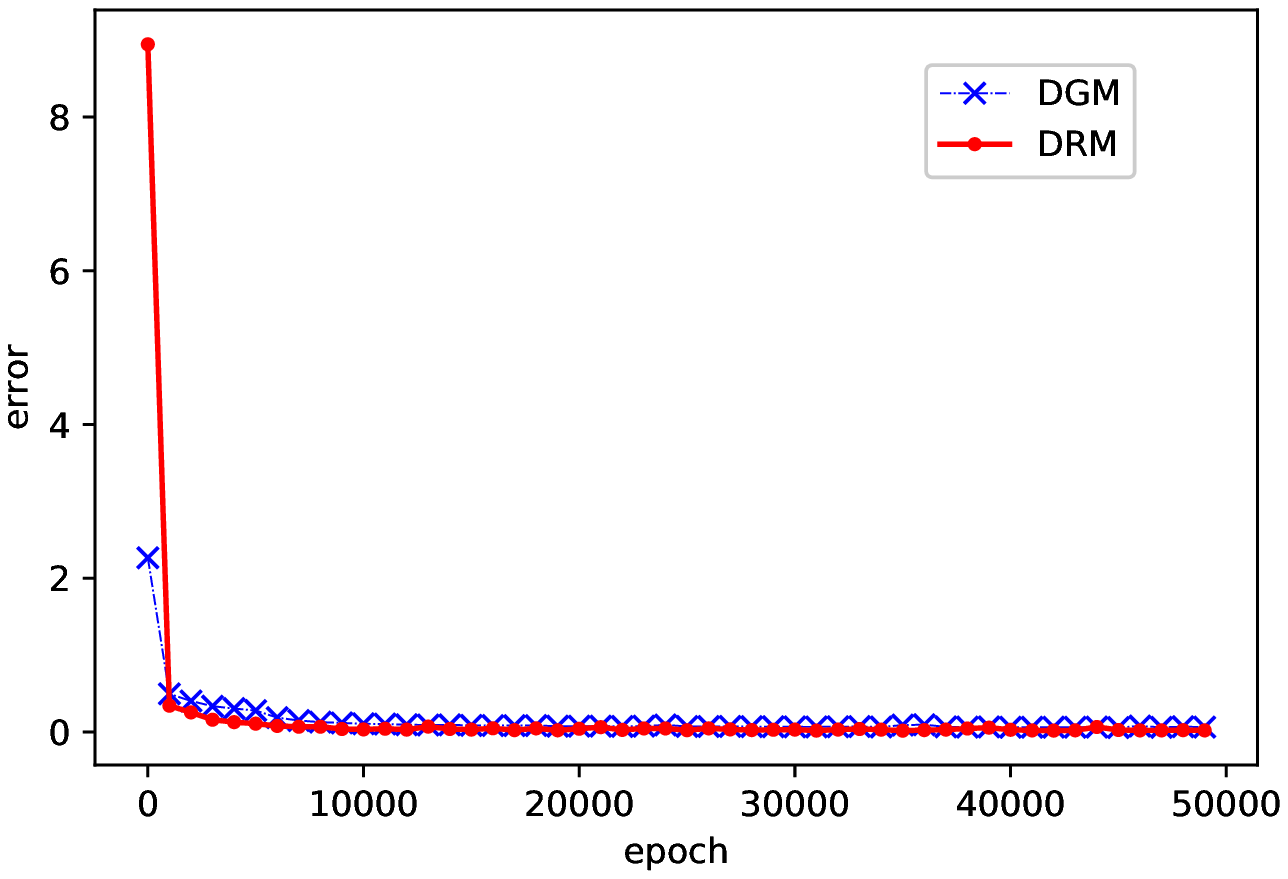}
		}
		\subfigure[Periodic]{
			\includegraphics[width=2.0in]{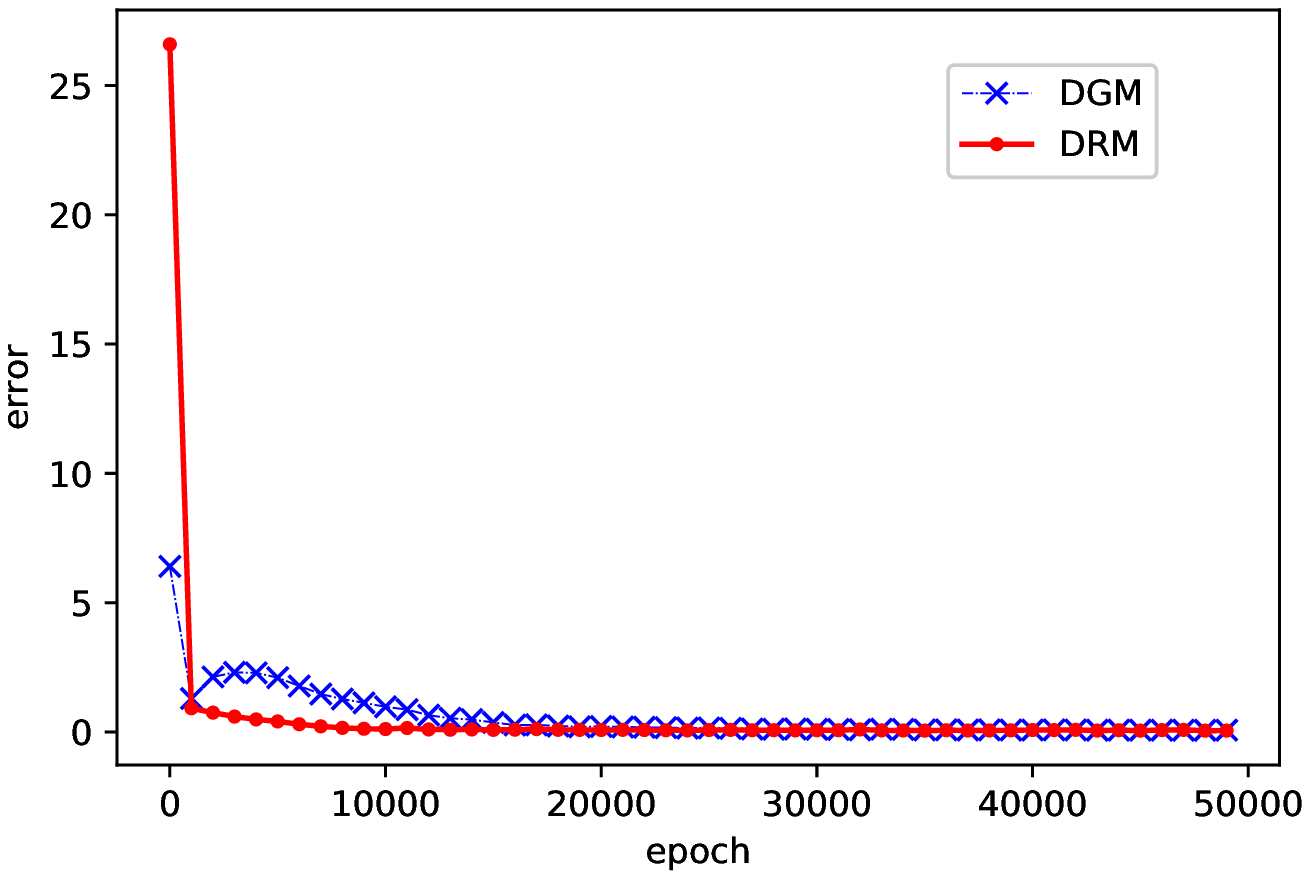}
		}
		\caption{Training processes of DGM and DRM for four boundary conditions in 8D. Each neural network contains three residual blocks with sixteen neural units in each layer. The mini-batch size is $2000$ in the domain for Dirichlet, Neumann, and Robin boundary conditions, and $4000$ in the domain for periodic boundary condition, $1600$ on the boundary for Dirichlet, Neumann, and Robin boundary conditions, and $16000$ on the boundary for periodic boundary condition. The penalty parameter $\lambda = 100.0$ for Dirichlet boundary condition, $\lambda = 1.0$ for Neumann boundary condition, $\lambda = 10.0$ for Robin boundary condition, and $\lambda_1 = 1.0, \lambda_2 = 0.5$ for periodic boundary condition.}
		\label{DGM v.s. DRM in 8D}
	\end{figure}

\begin{figure}[H]
		\centering
		\subfigure[Dirichlet]{
			\includegraphics[width=2.0in]{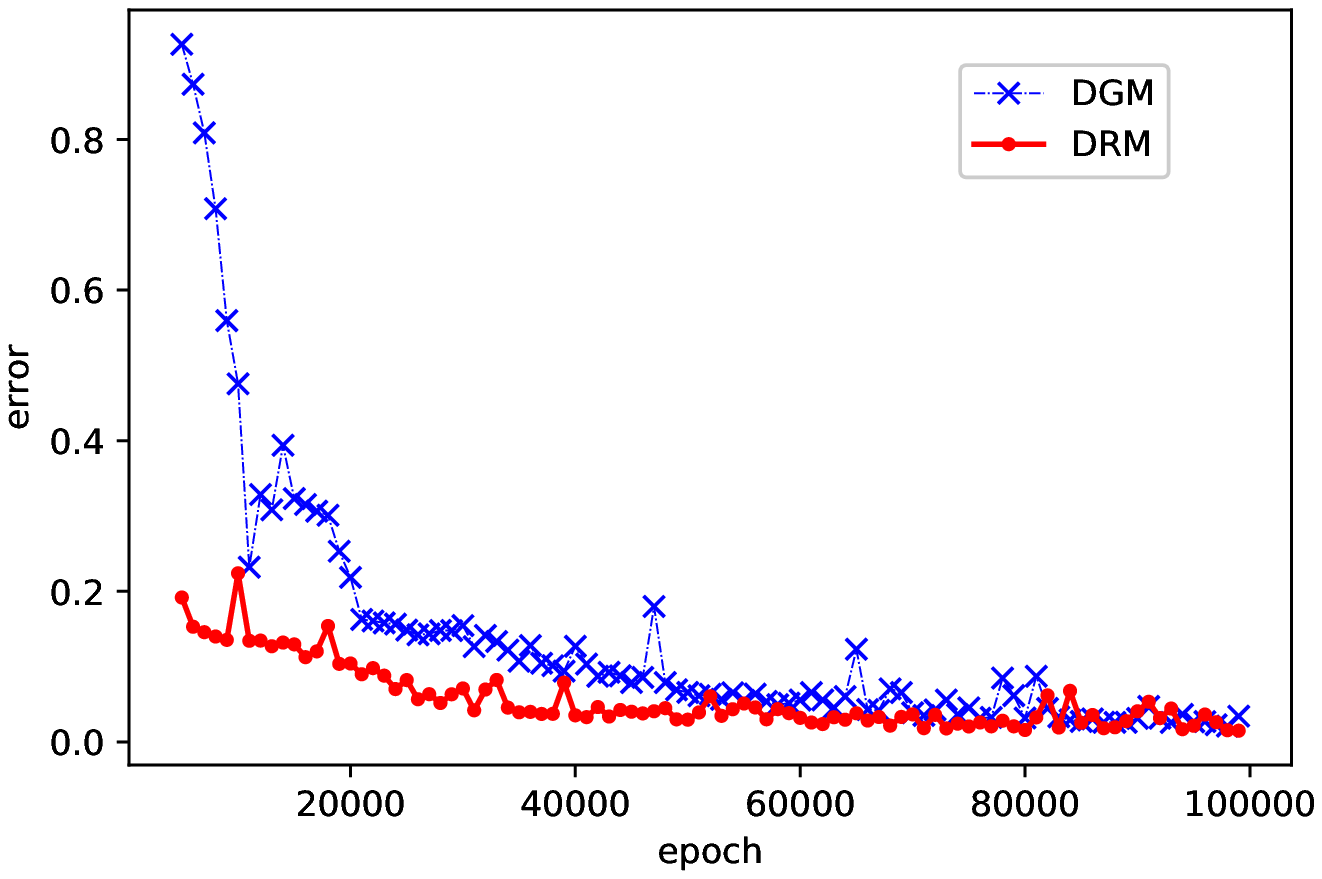}
		}
		\subfigure[Neumann]{
			\includegraphics[width=2.0in]{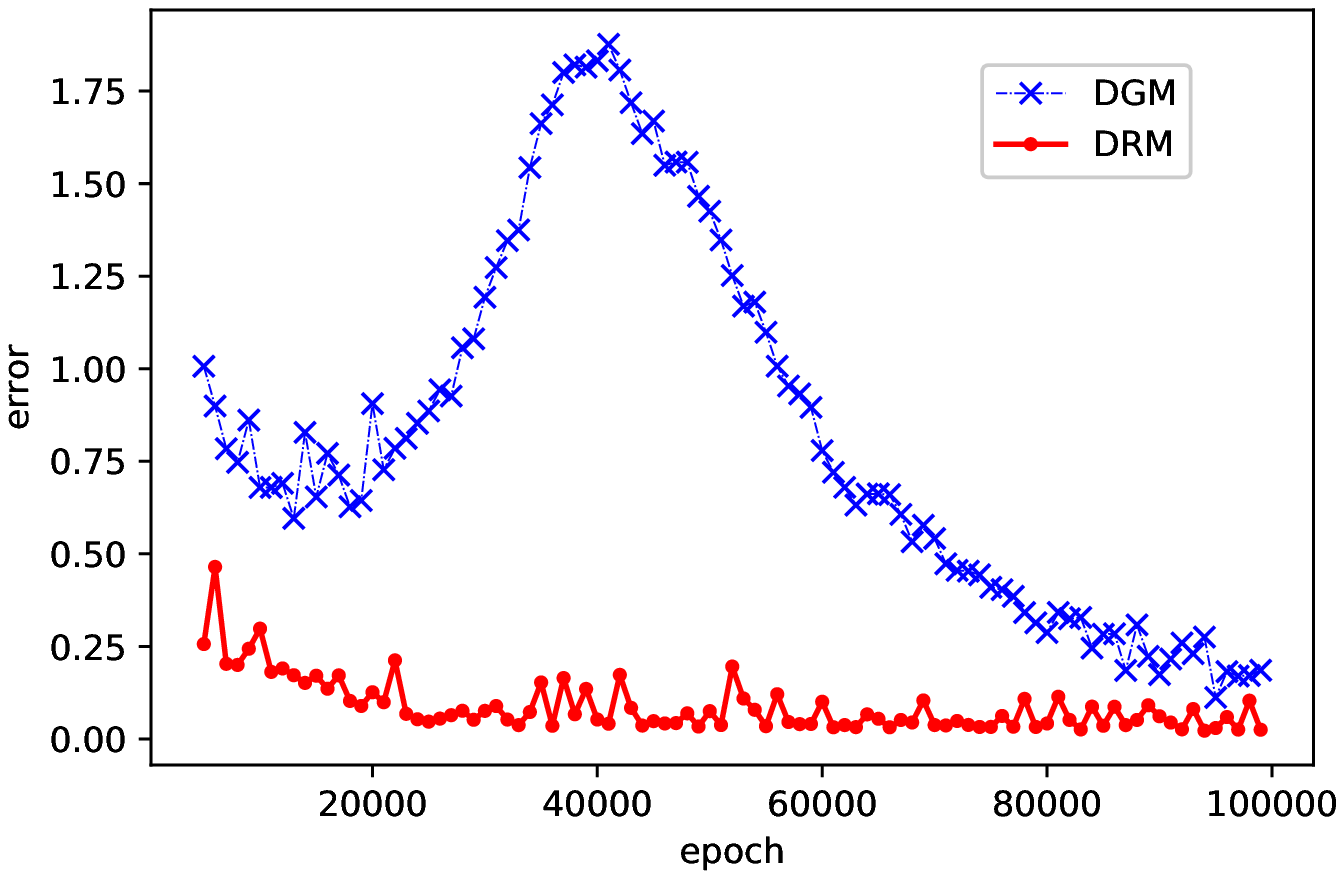}
		}
		\quad    
		\subfigure[Robin]{
			\includegraphics[width=2.0in]{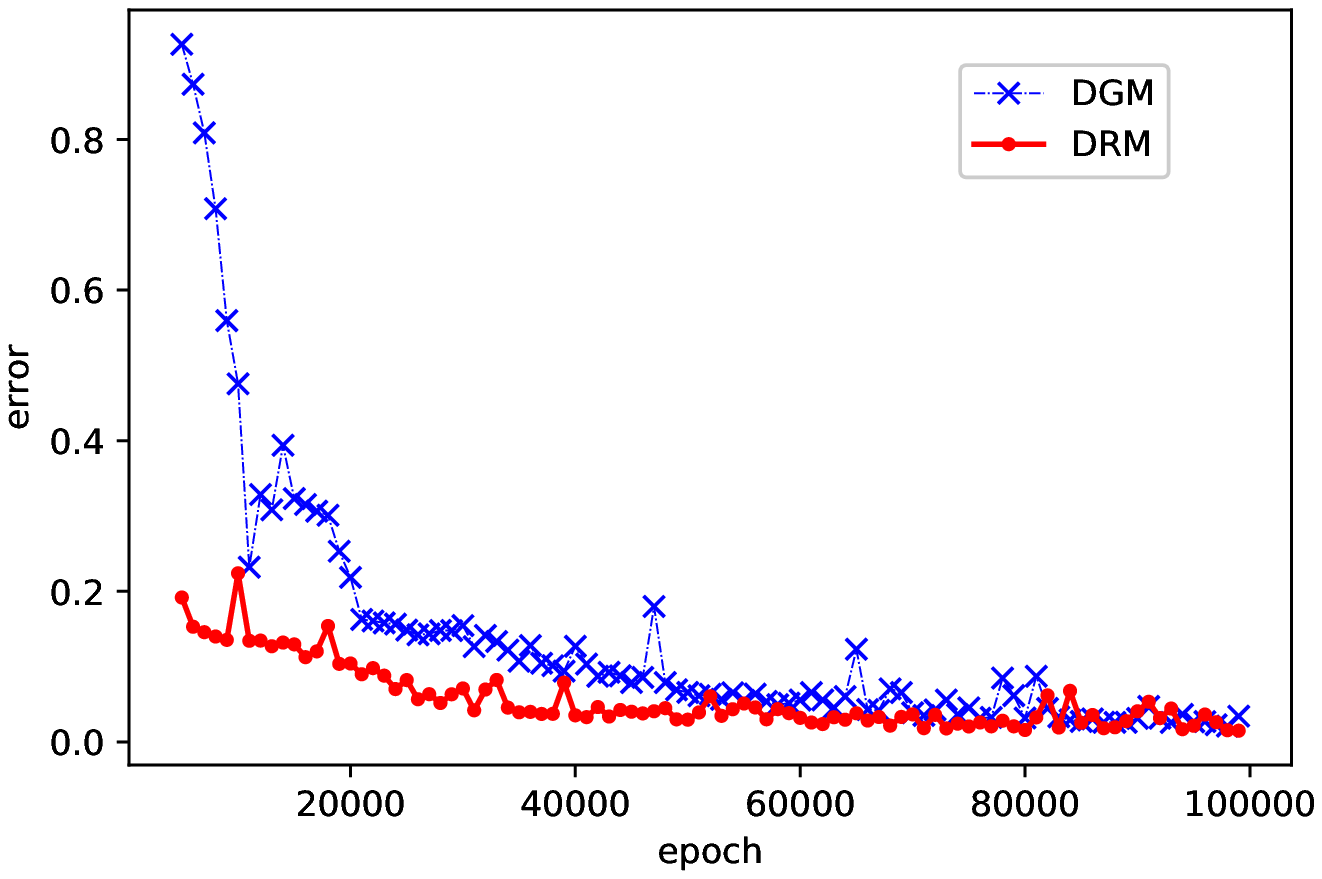}
		}
		\subfigure[Periodic]{
			\includegraphics[width=2.0in]{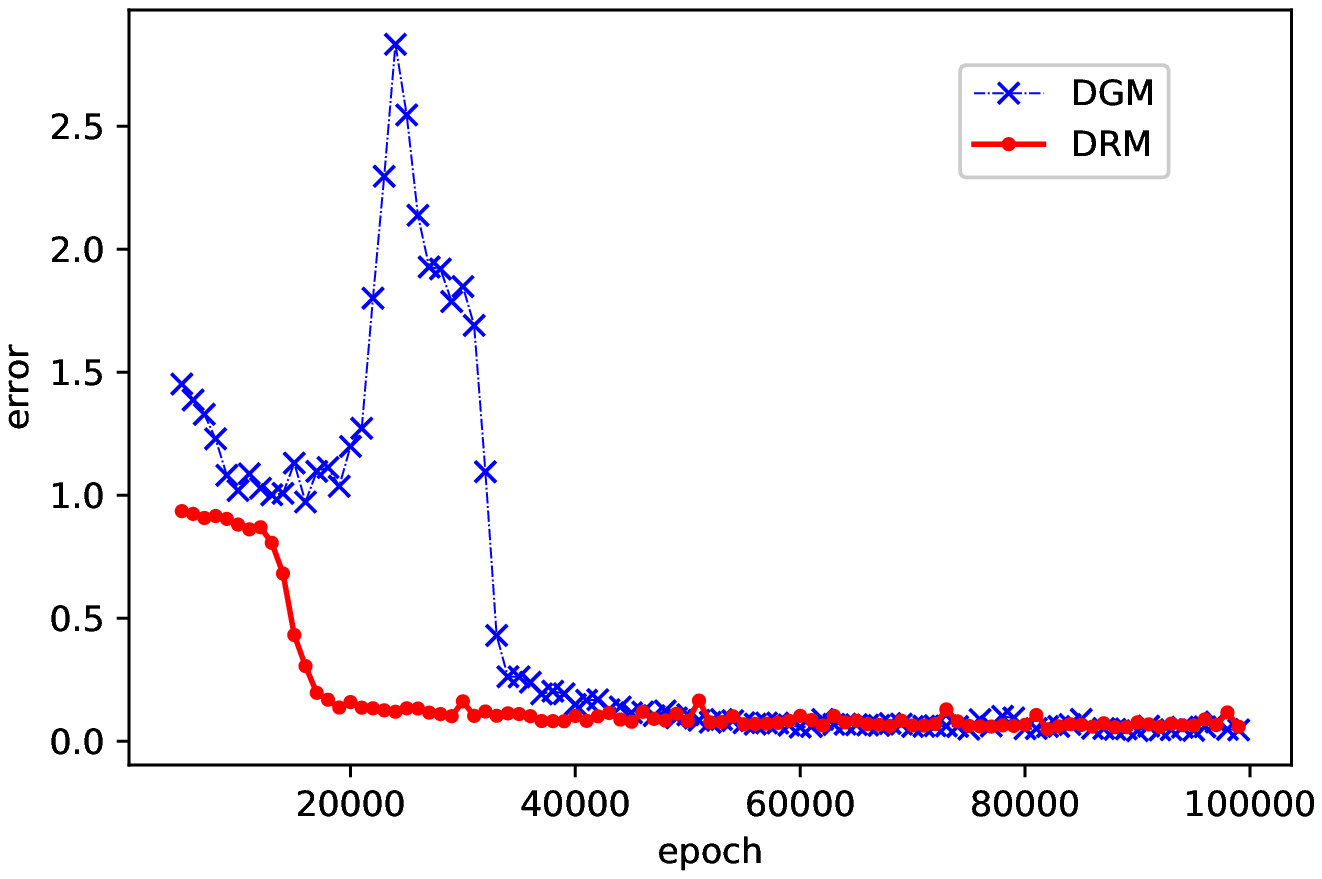}
		}
		\caption{Training processes of DGM and DRM for four boundary conditions in 16D. Each neural network contains three residual blocks with thirty-two neural units in each layer. The mini-batch size is $2000$ in the domain and $3200$ on the boundary. The penalty parameter $\lambda = 100.0$ for Dirichlet boundary condition, $\lambda = 1.0$ for Neumann boundary condition, $\lambda = 10.0$ for Robin boundary condition, and $\lambda_1 = 10.0, \lambda_2 = 5.0$ for periodic boundary condition.}
		\label{DGM v.s. DRM in 16D}
	\end{figure}

	Figure \ref{DGM v.s. DRM in 2D} - Figure \ref{DGM v.s. DRM in 16D} record the training processes of DGM and DRM in 2D, 4D, 8D, and 16D, respectively. One general trend we have observed is that DGM converges faster than DRM in the low-dimensional case; see 2D for example, while it is the opposite in the high-dimensional case; see 16D for example. In lower dimensions, both DGM and DRM converge well. However, in 16D, a significant amount of efforts have been paid in order to achieve the convergence in DGM. From Figure \ref{DGM v.s. DRM in 2D} - Figure \ref{DGM v.s. DRM in 16D}, a general observation is that the penalty parameter decreases from Dirichlet, Neumann, Robin, to periodic boundary conditions. Since this parameter is a bit tuned to get a better approximation for a given DNN, a larger damping parameter implies a better agreement between the DNN solution and the exact solution on the boundary.

	{Table \ref{error} records relative $L^2$ errors for DGM and DRM when $d = 2, 4, 8, 16, 32$. The maximum number of training epochs is set to be $10000$ in 2D, $20000$ in 4D, $50000$ in 8D, $100000$ in 16D, and $400000$ in 32D. Other parameters of neural network are the same as those in Figure \ref{DGM v.s. DRM in 2D} - Figure \ref{DGM v.s. DRM in 16D}. When $d=32$, the batch size is $2000$ in the domain and $3200$ on the boundary. The network width is $32$ and the depth is $4$. $-$ means the training process does not converge after the number of training epochs has reached the maximum epoch number.} Generally speaking, DGM has a better approximation accuracy in low-dimensional cases; see 2D and 4D for example, while DRM outperforms in high-dimensional cases; see 8D and 16D for example. For \eqref{eqn:poisson}, the second-order derivative appears in the formulation of DGM while only the first-order derivative exists in DRM. Therefore, to some extent, this is out of expectation since the exact solution here is smooth and DGM should approximate the exact solution better.
	\begin{table}[htbp]
		\centering
		\caption{Relative $L^2$ errors for four different boundary conditions in different dimensions. The number of training epochs is $10000$ in 2D, $20000$ in 4D, $50000$ in 8D, $100000$ in 16D, and $200000$ in 32D. Other parameters in DNNs are specified in Figure \ref{DGM v.s. DRM in 2D} - Figure \ref{DGM v.s. DRM in 16D}.}
		\label{error}      
		\begin{tabular}{c|cc|cc|cc|cc}
			\toprule[2pt]
			\noalign{\smallskip}
			\multirow{2}*{$d$} 
			&\multicolumn{2}{c}{Dirichlet} 
			&\multicolumn{2}{c}{Neumann}
			&\multicolumn{2}{c}{Robin} 
			&\multicolumn{2}{c}{Periodic} \\			
			&\multicolumn{1}{c}{DGM} & \multicolumn{1}{c}{DRM}
			&\multicolumn{1}{c}{DGM} & \multicolumn{1}{c}{DRM}
			&\multicolumn{1}{c}{DGM} & \multicolumn{1}{c}{DRM}
			&\multicolumn{1}{c}{DGM} & \multicolumn{1}{c}{DRM} \\
			\noalign{\smallskip}
			\midrule[1pt]
			\noalign{\smallskip}
			\multirow{1}*{2}
			& 0.0071 & 0.0236 & 0.0020 & 0.0078  & 0.0006 & 0.0065 & 0.0063 & 0.0115 \\
			\multirow{1}*{4} 
			& 0.0074 & 0.0105 & 0.0128 & 0.0336  & 0.0197 & 0.0622 & 0.0449 & 0.0514 \\
			\multirow{1}*{8} 
			& 0.0226 & 0.0256 & 0.0674 & 0.0199  & 0.0561 & 0.0221 & 0.0672 & 0.0573 \\
			\multirow{1}*{16} 
			& 0.0290 & 0.0224 & 0.1747 & 0.0368  & 0.0938 & 0.0379 & 0.0525 & 0.0617 \\
			\multirow{1}*{{32}} 
			& {0.0912} & {0.0561} & {-} & {0.0399}  & {0.1828} & {0.0303} & {-} & {-} \\			
			\noalign{\smallskip}
			\bottomrule[2pt]
		\end{tabular}
	\end{table}

	\subsection{Dependence on network structures}
	The above observations hold true over a wide range of issues, such as penalty parameter, mini-batch size, activation function, neural depth, and neural width. We will show how the approximation accuracy of DGM and DRM depends on these issues by several representative results in what follows.
	
	\subsubsection{Penalty parameter}
	Consider Dirichlet boundary condition in 4D. Relative $L^2$ errors of DGM and DRM are recorded in Table \ref{4D Dirichlet penalty} for different penalty parameters. In theory, the damping parameter $\lambda$ shall be infinity if the exact solution is found. In practice, instead, for a given DNN, $\lambda$ shall always be a finite number. It is observed from Table  \ref{4D Dirichlet penalty} that the larger the penalty parameter $\lambda$ is, the better the approximation is. However, if $\lambda$ is set to be too small or too large, then the penalty term can be ignored or be dominant{; see Table \ref{penalty 2D} for example where the approximation error increases if $\lambda$ is too large}. This may result in wrong DNN solutions, i.e., a DNN approximation satisfies the PDE but not the boundary condition or satisfies the boundary condition but not the PDE. Therefore, for a given DNN, how to choose a penalty parameter which grantees the optimal approximation accuracy is of particular importance and deserves further consideration.
	\begin{table}[H]
		\centering 
		\caption{Relative $L^2$ errors of DGM and DRM in terms of penalty parameter $\lambda$ for Dirichlet boundary condition in 4D. The neural network contains three residual blocks with eight neural units in each layer. The activation function is $swish(x)$. The mini-batch size is $2000$ in the domain and $800$ on the boundary.}
		\begin{tabular}{c|c|c}
			\toprule[2pt]
			$\lambda$ & DGM & DRM \\ 
			\toprule[2pt]
			0.1               & 0.2186       & 0.0185       \\
			1.0               & 0.0366       & 0.0176       \\
			10.0              & 0.0127       & 0.0196       \\
			100.0             & 0.0081       & 0.0083       \\ 
			\toprule[2pt]
		\end{tabular}
		\label{4D Dirichlet penalty}
	\end{table}

	\subsubsection{Mini-batch size}
	Consider Robin boundary condition in 4D. {Relative $L^2$ errors of DGM and DRM are recorded in Table \ref{4D Robin mini-batch} for different mini-batch sizes in the domain with fixed mini-batch size on the boundary and  in Table \ref{4D Robin mini-batch  boundary} for different mini-batch size on the boundary and fixed mini-batch size in the domain, respectively. From the results, one can expect that a balance of sampling points in the domain and on the boundary will yield the optimal approximation accuracy. Interested readers may refer to \cite{2020arXiv200206269V} for such an effort.} 
	\begin{table}[H]
		\centering
		\caption{Relative $L^2$ errors of DGM and DRM in terms of mini-batch size in the domain. The neural network contains three residual blocks with eight neural units in each layer. The activation function is $swish(x)$. The mini-batch size is $800$ on the boundary. The penalty parameter $\lambda = 100$.}
		\begin{tabular}{c|c|c}
			\toprule[2pt]
			Mini-batch size & DGM & DRM \\ 
			\toprule[2pt]
			500             & 0.0822       & 0.0230       \\
			1000            & 0.1064       & 0.0266       \\
			2000            & 0.0197       & 0.0622       \\
			4000            & 0.1026       & 0.0321       \\ 
			\toprule[2pt]
		\end{tabular}
		\label{4D Robin mini-batch} 
	\end{table}
	
	\begin{table}[H]
		\centering
		\caption{{Relative $L^2$ errors of DGM and DRM in terms of mini-batch size on the boundary. The neural network contains three residual blocks with eight neural units in each layer. The activation function is $swish(x)$. The mini-batch size is $500$ in the domain. The penalty parameter $\lambda = 100$.}}
		\begin{tabular}{c|c|c}
			\toprule[2pt]
			{Mini-batch size} & {DGM} & {DRM} \\ 
			\toprule[2pt]
			{400}            & {0.0339}       & {0.0184}       \\
			{800}             & {0.1048}       & {0.0204}       \\
			{1600}            & {0.5105}       & {0.0409}       \\
			\toprule[2pt]
		\end{tabular}
		\label{4D Robin mini-batch boundary} 
	\end{table}

	\subsubsection{Activation function}
	Consider Neumann boundary condition in 4D. Table \ref{4D Neumann acti} records relative $L^2$ errors of DGM and DRM in terms of several activation functions. From Table \ref{4D Neumann acti}, it is recognized that the choice of activation functions is quite important. The failure of $relu(x)$ in DGM is due to the low regularity of the activation function and the usage of higher derivatives in the loss function of DGM. This is why a better performance of DGM for smooth solutions is expected while DRM is expected to be better for low-regularity solutions. Based on results in Section \ref{sec:bcresult}, we know that the former is not always true. In Section \ref{sec:nonlinearresult}, the latter is not true as well.
	\begin{table}[H]
		\centering
		\caption{Relative $L^2$ errors in terms of activation function. The neural network contains three residual blocks with eight neural units in each layer. The mini-batch size is $2000$ in the domain and $800$ on the boundary. The penalty parameter $\lambda = 500$.}
		\begin{tabular}{c|c|c}
			\toprule[2pt]
			Activation function & DGM & DRM \\ 
			\toprule[2pt]
			$relu(x)$           & 0.9992       & 0.0783       \\
			$sigmoid(x)$        & 0.0226       & 0.0136       \\
			$swish(x)$          & 0.0176       & 0.0169       \\
			${(\sin x)}^3$      & 0.0231       & 0.0110       \\
			{$swish(ax)$}          & {0.0147}       & {0.0112}       \\ 
			\toprule[2pt]
		\end{tabular}
		\label{4D Neumann acti}
	\end{table}
	
	\subsubsection{Neural depth and neural width}
	Consider Dirichlet boundary condition in 4D. Table \ref{4D Dirichlet depth} and Table \ref{4D Dirichlet width} record relative $L^2$ errors of DGM and DRM in terms of neural depth $n$ and neural width $m$, respectively. It is expected that approximation errors of DGM and DRM reduce as $n$ and $m$ increase to some extent. Unlike classical numerical methods, a systematic reduction of errors cannot be observed for DNNs.
	\begin{table}[H]
		\centering
		\caption{Relative $L^2$ errors in terms of neural depth $n$. Each neural network contains varying residual blocks with eight neural units in each layer. The activation function is $swish(x)$. The mini-batch size is $2000$ in the domain and $800$ on the boundary. The penalty parameter $\lambda = 100$.}
		\begin{tabular}{c|c|c}
			\toprule[2pt]
			Neural depth $n$ & DGM & DRM \\ 
			\toprule[2pt]
			2                         & 0.0114       & 0.0193       \\
			3                         & 0.0074       & 0.0105       \\
			4                         & 0.0108       & 0.0057       \\
			\toprule[2pt]
		\end{tabular}
		\label{4D Dirichlet depth}
	\end{table}
	\begin{table}[H]
		\centering
		\caption{Relative $L^2$ errors in terms of neural width $m$. Each neural network contains three residual blocks with varying neural units in each layer. The activation function is $swish(x)$. The mini-batch size is $2000$ in the domain and $800$ on the boundary. The penalty parameter $\lambda = 100$.}
		\begin{tabular}{c|c|c}
			\toprule[2pt]
			Neural width $m$             & DGM & DRM \\ 
			\toprule[2pt]
			4                         & 0.0218       & 0.1118       \\
			6                         & 0.0208       & 0.1124       \\
			8                         & 0.0074       & 0.0105       \\
			10                        & 0.0072       & 0.0095       \\
			\toprule[2pt]
		\end{tabular}
		\label{4D Dirichlet width}
	\end{table}

	\subsubsection{Periodic boundary condition}\label{sec:pbc}
	
	{Consider equation \eqref{eqn:poisson} satisfying periodic boundary condition with period $p_i=2, i=1,\cdots,d$. In order to construct a DNN which satisfies the periodicity, following \cite{han2020solving}, we construct a transform $\mathbb{R}^d \to  \mathbb{R}^{2kd}$ for the input $x = (x_1,\cdots,x_d)$ before the first fully connected layer of the neural network. The component $x_i$ of $x$ is transformed as follows
	\begin{equation*}
		x_i  \to \{\sin(j \cdot 2\pi \frac {x_i} { p_i}) , \cos(j \cdot 2\pi \frac {x_i} { p_i})\}_{j = 1}^k
	\end{equation*}
	for $i = 1, \cdots, d$.}
	
	{Since the exact solution in Section \ref{sec:pbcpenalty} can be exactly expressed by the above transform and the approximation error is significantly small. To avoid this, we choose the exact solution $u(x) = \sum_{i = 1}^d \cos(\pi x_i) \cos(2 \pi x_i) $, which cannot be explicitly represented by the above transform. The approximation error is recorded in Table \ref{p1error}.}
	\begin{table}[H]
		\centering 
		\caption{{Relative $L^2$  errors for periodic boundary condition without the penalty term in different dimensions. Here we set $k = 3$.}}
		\begin{tabular}{|c|c|c|c|}
			\hline
			{d} & {DGM} & {DRM} \\
			\hline
			{2}  & {0.0033}  & {0.0281} \\
			\hline
			{4}  & {0.0012}  & {0.0656} \\
			\hline
			{8}  & {0.0021}  & {0.0657} \\
			\hline
			{16} & {0.0067}  & {0.0490} \\
			\hline
		\end{tabular}
		\label{p1error}
	\end{table}

	\subsection{A nonlinear problem with low-regularity solution}\label{sec:nonlinearresult}
	
	Note that all the previous examples are linear PDEs and their solutions belong to $C^{\infty}(\Omega)$. Next, we study a nonlinear PDE with the low-regularity solution. The nonlinear problem over the unit sphere $\Omega = \{x \in \mathbb{R}^d:|x| < 1\} $ reads as
	\begin{equation}
	\begin{cases}
	- \Delta u +  u^3 = f(x), & \; \text{in} \; \Omega,\\
	u(x) = 0, & \; \text{on} \; \partial \Omega.
	\end{cases}
	\end{equation}
	The exact solution $u(x) = \sin \left(\frac{\pi}{2} (1 - |x|)\right)  \in C^1(\Omega)$ but $u(x) \notin C^2(\Omega)$, and
	\begin{equation*}
		f(x) = \frac{\pi^2}{4} \sin \left(\frac{\pi}{2} (1 - |x|)\right) + \frac{\pi}{2} \cos \left(\frac{\pi}{2} (1 - |x|)\right) \frac{d-1}{|x|} +  \sin^3 \left(\frac{\pi}{2} (1 - |x|)\right).
	\end{equation*}
	Loss functions associated to DGM and DRM are
	\begin{equation}
	\mathcal{J}_{\textrm{DGM}} [u(x;\theta)] = \int_{\Omega} {| - \Delta u(x;\theta) +  {u(x;\theta)}^3 - f(x) |}^2 \mathrm{d}x,
	\end{equation}
	\begin{equation}
	\mathcal{J}_{\textrm{DRM}} [u(x;\theta)] = \int_{\Omega} {\frac{1}{2}  {|\nabla u(x;\theta)|}^2 + \frac{1}{4} {u(x;\theta)}^4  - f(x) u(x;\theta)} \mathrm{d}x,
	\end{equation}
	and the penalty term is
	\begin{equation}
	\mathcal{B}_{\textrm{D}} [u(x;\theta)] = \int_{\partial\Omega} {| u(x;\theta)|}^2 \mathrm{d}s.
	\end{equation}
	Thus, total loss functions of DGM and DRM with penalty are  
	\begin{equation}
		\mathcal{I}_{\textrm{DGM}} [u(x;\theta)]  = \mathcal{J}_{\textrm{DGM}} [u(x;\theta)] + \lambda \mathcal{B}_{\textrm{D}} [u(x;\theta)],
	\end{equation}
	and
	\begin{equation}
		\mathcal{I}_{\textrm{DRM}} [u(x;\theta)]  = \mathcal{J}_{\textrm{DRM}} [u(x;\theta)] + \lambda \mathcal{B}_{\textrm{D}} [u(x;\theta)],
	\end{equation}
	respectively.

	\subsubsection{Dimensional dependence}
	Relative $L^2$ errors of DGM and DRM are reported in different dimensions $d = 2, 4, 8$. Each neural network contains three residual blocks with varying neural units in each layer. The number of neural units is 8 for 2D and 4D, and 16 for 8D. The activation function is $swish(x)$. The mini-batch size is $2000$ in the domain and $400$ on the boundary in 2D, $1000$ in the domain and $800$ on the boundary in 4D, and $1000$ in the domain and $1600$ on the boundary for 8D. The penalty parameter is $50.0$ for 2D, $100.0$ for 4D, and $400.0$ for 8D. To our surprise, DGM outperforms DRM by over one order of magnitude. Note that the exact solution is only in $C^1(\Omega)$, the second-order derivative in space appears in DGM while only the first-order derivative in space is needed in DRM. Therefore, such an observation definitely deserves further investigation. Moreover, this observation holds true over a wide range of issues, such as penalty parameter, mini-batch size, activation function, neural depth, and neural width. We will show how the approximation accuracy of DGM and DRM depends on a couple of representative issues in what follows.
	\begin{table}[H]
		\centering 
		\caption{Relative $L^2$ errors of DGM and DRM in different dimensions. Each neural network contains three residual blocks with varying neural units in each layer. The number of neural units is 8 for 2D and 4D, and 16 for 8D. The activation function is $swish(x)$. The mini-batch size is $2000$ in the domain and $400$ on the boundary in 2D, $1000$ in the domain and $800$ on the boundary in 4D, and $1000$ in the domain and $1600$ on the boundary for 8D. The penalty parameter is $50.0$ for 2D, $100.0$ for 4D, and $400.0$ for 8D.}
		\begin{tabular}{c|c|c}
			\toprule[2pt]
			$d$      & DGM & DRM \\ 
			\toprule[2pt]
			2              & 0.0003       & 0.0090       \\
			4              & 0.0055       & 0.0777       \\
			8              & 0.0292       & 0.1603       \\
			\toprule[2pt]
		\end{tabular}
		\label{Dirichlet Dim 3}
	\end{table}
	
	\subsubsection{Penalty parameter}
	Table \ref{penalty 2D} records relative $L^2$ errors of DGM and DRM in terms of penalty parameter $\lambda$ in 2D. 
	
	\begin{table}[H]
		\centering 
		\caption{Relative $L^2$ errors of DGM and DRM in terms of penalty parameter $\lambda$ in 2D. Each neural network contains three residual blocks with eight neural units in each layer. The activation function is $swish(x)$. The mini-batch size is $2000$ in the domain and $400$ on the boundary.}
		\begin{tabular}{c|c|c}
			\toprule[2pt]
			$\lambda$ & DGM & DRM \\ 
			\toprule[2pt]
			50.0              & 0.0011       & 0.0517       \\
			100.0             & 0.0022       & 0.0161       \\
			200.0             & 0.0015       & 0.0076       \\
			400.0             & 0.0003       & 0.0090       \\ 
			{2000.0}           & {0.0006}       & {0.0175}       \\ 
			{10000.0}           & {0.0038}       & {0.0300}       \\
			{100000.0}          & {0.0106}        & {0.3873}        \\
			\toprule[2pt]
		\end{tabular}
		\label{penalty 2D}
	\end{table}

	\subsubsection{Activation function}
	Table \ref{4D acti} records relative $L^2$ errors of DGM and DRM with respect to activation function in 4D. From Table \ref{4D acti}, we see that $relu(x)$ still has some problem due to the same reason and $swish(x)$ is the best among all the tested functions.
	\begin{table}[H]
		\centering
		\caption{Relative $L^2$ errors of DGM and DRM with respect to activation function in 4D. Each neural network contains three residual blocks with eight neural units in each layer. The mini-batch size is $1000$ in the domain and $800$ on the boundary. The penalty parameter $\lambda$ is $100$.}
		\begin{tabular}{c|c|c}
			\toprule[2pt]
			Activation function & DGM & DRM \\ 
			\toprule[2pt]
			$relu(x)$           & 0.9990       & 0.1546       \\
			$sigmoid(x)$        & 0.0262       & 0.0881       \\
			$swish(x)$          & 0.0055       & 0.0777       \\
			${(\sin x)}^3$      & 0.0146       & 0.0907       \\ 
			\toprule[2pt]
		\end{tabular}
		\label{4D acti}
	\end{table}

	\subsubsection{With versus without penalty}
	For Dirichlet boundary condition, as discussed earlier, we can actually avoid the penalty term \cite{berg2018unified} by constructing a trail function in the  form of \eqref{eqn:twostage}. Since $\Omega$ is a unit sphere, there exists a simple way to construct a trail function which automatically satisfies the exact boundary condition. Precisely, we can build the neural network solution in the form of $u(x;\theta) = (1 - |x|)  DNN(x;\theta)$, where $DNN(x;\theta)$ is the DNN approximation to be trained. This will be used for both DGM and DRM without penalty term for the comparison purpose.
	
\begin{figure}[H]
		\centering
		\includegraphics[width= 12cm]{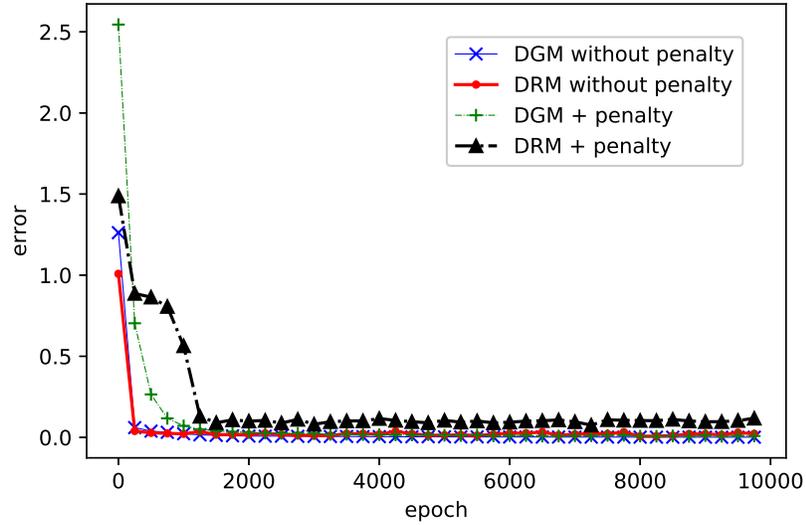}
		\caption{Training processes of DGM and DRM with or without penalty. Each neural network contains three residual blocks with eight neural units in each layer. The mini-batch size is $1000$ in the domain and $800$ on the boundary. The penalty parameter $\lambda$ is $100.0$. The total number of epochs is $10000$.}
		\label{nonlinearPDE}
	\end{figure}

	Figure \ref{nonlinearPDE} plots training processes of DGM and DRM with or without penalty and Table \ref{nonlinear Table} records the corresponding relative $L^2$ errors in 4D. Each neural network contains three residual blocks with eight neural units in each layer. The mini-batch size is $1000$ in the domain and $800$ on the boundary. The penalty parameter $\lambda$ is $100.0$. From Figure \ref{nonlinearPDE}, without penalty, we see that both DGM and DRM converge better. Sometimes we even see that DGM and DRM without penalty converge while do not converge in the presence of penalty term. Besides, from Table \ref{nonlinear Table}, we see that DGM outperforms DRM by over one order of magnitude regardless of the penalty term, and both methods perform better by over one order of magnitude if the trail function automatically satisfies the boundary condition. These together show the great importance of boundary conditions. A better treatment not only facilitates the training process but also provides a better approximation accuracy for the same network setup.
	\begin{table}[H]
		\centering 
		\caption{Relative $L^2$ erros of DGM and DRM with or without penalty.}
		\begin{tabular}{c|c|c}
			\toprule[2pt]
			With or without penalty & DGM & DRM \\ 
			\toprule[2pt]
			With penalty        & 0.0055       & 0.0777       \\
			Without penalty    & 0.0002       & 0.0084       \\
			\toprule[2pt]
		\end{tabular}
		\label{nonlinear Table}
	\end{table}
	
	\section{Conclusions}\label{sec:conclusion}

	In the work, we have conducted a comprehensive study of four different boundary conditions, i.e., Dirichlet, Neumann, Robin, and periodic boundary conditions, using two representative methods: DRM and DGM. It is thought that DGM works better for smooth solutions while DRM works better for low-regularity solutions. However, by a number of examples, we have observed that DRM can outperform DGM with a clear dependence of dimensionality even for smooth solutions and DGM can also outperform DRM for low-regularity solutions. Besides, in some cases, when the boundary condition can be implemented in an exact manner, we have found that such a strategy not only provides a better approximate solution but also facilitates the training process.
	
    There are several interesting issues which deserves further considerations. Since the penalty method works in general, the most important one is the choice of penalty parameters. For a fixed neural structure, a good choice of these parameters not only facilitates the training process but also provides a better approximation. Another issue is to understand why DGM outperforms DRM for low-regularity problems.
	
	\section*{Acknowledgements}
	This work is supported in part by the grants NSFC 11971021 and National Key R\&D Program of China (No. 2018YF645B0204404) (J.~Chen), NSFC 11501399 (R.~Du). We are grateful to Liyao Lyu for helpful discussions. All codes for producing the results in this work are available at \url{https://github.com/wukekever/DGM-and-DRM}.

	\newpage
	\bibliographystyle{unsrt}
	\bibliography{refs}

\begin{thebibliography}{10}

\bibitem{NIPS2012_4824}
Alex Krizhevsky, Ilya Sutskever, and Geoffrey~E. Hinton.
\newblock Imagenet classification with deep convolutional neural networks.
\newblock {\em Communications of the ACM}, 60:84--90, 2012.

\bibitem{hinton2012deep}
Hinton Geoffrey, Deng Li, Yu~Dong, Dahl George, Mohamed Abdel-rahman, Jaitly
  Navdeep, Senior Andrew, Vanhoucke Vincent, Nguyen Patrick, Kingsbury Brian,
  and Sainath Tara.
\newblock Deep neural networks for acoustic modeling in speech recognition.
\newblock {\em IEEE Signal Processing Magazine}, 29:82--97, 2012.

\bibitem{Goodfellow2016}
Ian Goodfellow, Yoshua Bengio, and Aaron Courville.
\newblock {\em Deep Learning}.
\newblock MIT Press, 2016.

\bibitem{lagaris1998artificial}
Isaac~E Lagaris, Aristidis Likas, and Dimitrios~I Fotiadis.
\newblock Artificial neural networks for solving ordinary and partial
  differential equations.
\newblock {\em IEEE transactions on neural networks}, 9(5):987--1000, 1998.

\bibitem{raissi2018hidden}
Maziar Raissi and George~Em Karniadakis.
\newblock Hidden physics models: Machine learning of nonlinear partial
  differential equations.
\newblock {\em Journal of Computational Physics}, 357:125--141, 2018.

\bibitem{weinan2017deep}
E~Weinan, Jiequn Han, and Arnulf Jentzen.
\newblock Deep learning-based numerical methods for high-dimensional parabolic
  partial differential equations and backward stochastic differential
  equations.
\newblock {\em Communications in Mathematics and Statistics}, 5(4):349--380,
  2017.

\bibitem{berg2018unified}
Jens Berg and Kaj Nystr{\"o}m.
\newblock A unified deep artificial neural network approach to partial
  differential equations in complex geometries.
\newblock {\em Neurocomputing}, 317:28--41, 2018.

\bibitem{E2018}
Weinan E and Bing Yu.
\newblock The deep ritz method: A deep learning-based numerical algorithm for
  solving variational problems.
\newblock {\em Communications in Mathematics and Statistics}, 6(1):1--12, 2018.

\bibitem{long2017pde}
Zichao Long, Yiping Lu, Xianzhong Ma, and Bin Dong.
\newblock Pde-net: Learning pdes from data.
\newblock {\em arXiv preprint arXiv:1710.09668}, 2017.

\bibitem{han2018solving}
Jiequn Han, Arnulf Jentzen, and E~Weinan.
\newblock Solving high-dimensional partial differential equations using deep
  learning.
\newblock {\em Proceedings of the National Academy of Sciences},
  115(34):8505--8510, 2018.

\bibitem{deepGalerkin2018}
Justin Sirignano and Konstantinos Spiliopoulos.
\newblock {DGM: A} deep learning algorithm for solving partial differential
  equations.
\newblock {\em Journal of Computational Physics}, 375:1339--1364, 2018.

\bibitem{zang2020weak}
Yaohua Zang, Gang Bao, Xiaojing Ye, and Haomin Zhou.
\newblock Weak adversarial networks for high-dimensional partial differential
  equations.
\newblock {\em Journal of Computational Physics}, page 109409, 2020.

\bibitem{leveque2007finite}
Randall~J LeVeque.
\newblock {\em Finite difference methods for ordinary and partial differential
  equations: steady-state and time-dependent problems}, volume~98.
\newblock Siam, 2007.

\bibitem{brenner2007mathematical}
Susanne Brenner and Ridgway Scott.
\newblock {\em The mathematical theory of finite element methods}, volume~15.
\newblock Springer Science \& Business Media, 2007.

\bibitem{kharazmi2019variational}
EHSAN Kharazmi, Zhongqiang Zhang, and GE~Karniadakis.
\newblock Variational physics-informed neural networks for solving partial
  differential equations.
\newblock {\em arXiv preprint arXiv:1912.00873}, 2019.

\bibitem{liao2019deep}
Yulei Liao and Pingbing Ming.
\newblock Deep nitsche method: Deep ritz method with essential boundary
  conditions.
\newblock {\em arXiv preprint arXiv:1912.01309}, 2019.

\bibitem{hutzenthaler2020proof}
Martin Hutzenthaler, Arnulf Jentzen, Thomas Kruse, and Tuan~Anh Nguyen.
\newblock A proof that rectified deep neural networks overcome the curse of
  dimensionality in the numerical approximation of semilinear heat equations.
\newblock {\em SN Partial Differential Equations and Applications}, 1:1--34,
  2020.

\bibitem{evans_2010}
Evans~C. Lawrence.
\newblock {\em Partial differential equations (second edition)}.
\newblock American Mathematical Society, 2010.

\bibitem{hornik1989multilayer}
Kurt Hornik, Maxwell Stinchcombe, Halbert White, et~al.
\newblock Multilayer feedforward networks are universal approximators.
\newblock {\em Neural networks}, 2(5):359--366, 1989.

\bibitem{DBLP:journals/corr/HeZRS15}
Kaiming He, Xiangyu Zhang, Shaoqing Ren, and Jian Sun.
\newblock Deep residual learning for image recognition.
\newblock {\em 2016 IEEE Conference on Computer Vision and Pattern Recognition
  (CVPR)}, 2:770--778, 2016.

\bibitem{JAGTAP2020109136}
Ameya~D. Jagtap, Kenji Kawaguchi, and George~Em Karniadakis.
\newblock Adaptive activation functions accelerate convergence in deep and
  physics-informed neural networks.
\newblock {\em Journal of Computational Physics}, 404:109136, 2020.

\bibitem{2019arXiv190912228J}
Ameya~D. {Jagtap}, Kenji {Kawaguchi}, and George~Em {Karniadakis}.
\newblock {Locally adaptive activation functions with slope recovery term for
  deep and physics-informed neural networks}.
\newblock {\em arXiv e-prints}, page arXiv:1909.12228, 2019.

\bibitem{chen2019quasi}
Jingrun Chen, Rui Du, Panchi Li, and Liyao Lyu.
\newblock Quasi-monte carlo sampling for machine-learning partial differential
  equations.
\newblock {\em arXiv preprint arXiv:1911.01612}, 2019.

\bibitem{Ogata1989}
Yosihiko Ogata.
\newblock A {Monte Carlo} method for high dimensional integration.
\newblock {\em Numerische Mathematik}, 55(2):137--157, 1989.

\bibitem{ADAM}
Diederik~P. Kingma and Jimmy Ba.
\newblock Adam: A method for stochastic optimization.
\newblock {\em CoRR}, 1412.6980, 2014.

\bibitem{sheng2020pfnn}
Hailong Sheng and Chao Yang.
\newblock {PFNN: A} penalty-free neural network method for solving a class of
  second-order boundary-value problems on complex geometries.
\newblock {\em arXiv preprint arXiv:2004.06490}, 2020.

\bibitem{han2020solving}
Jiequn Han, Jianfeng Lu, and Mo~Zhou.
\newblock Solving high-dimensional eigenvalue problems using deep neural
  networks: A diffusion monte carlo like approach.
\newblock {\em arXiv preprint arXiv:2002.02600}, 2020.

\bibitem{2020arXiv200206269V}
Remco {van der Meer}, Cornelis {Oosterlee}, and Anastasia {Borovykh}.
\newblock {Optimally weighted loss functions for solving PDEs with Neural
  Networks}.
\newblock {\em arXiv e-prints}, page arXiv:2002.06269, 2020.

\end{thebibliography}

\end{document}